\documentclass[graybox]{svmult}
\pdfoutput=1

\usepackage{type1cm}        
\usepackage{makeidx}         
\usepackage{graphicx}        
\usepackage{multicol}        
\usepackage[bottom]{footmisc}

\usepackage{newtxtext}       %
\usepackage{newtxmath}       

\usepackage{combelow}

\usepackage{amsmath}
\usepackage{bbm}
\usepackage{euscript}

\def\mydate{\number\day\ {\ifcase\month \or January\or February\or
                   March\or April\or May\or June\or July\or August\or
                   September\or October\or November\or December\fi} \number\year}

\mathcode`\@="8000
{\catcode`\@=\active \gdef@{\mkern1mu}}     

 \renewcommand\mathcal{\EuScript}

\makeatletter
\renewcommand{\eqnarray}{%
   \stepcounter{equation}%
   \def\@currentlabel{\p@equation\theequation}%
   \global\@eqnswtrue
   \m@th
   \global\@eqcnt\z@
   \tabskip\@centering
   \let\\\@eqncr
   $$\everycr{}\halign to\displaywidth\bgroup
       \hskip\@centering$\displaystyle\tabskip\z@skip{##}$\@eqnsel
      &\global\@eqcnt\@ne\hskip0.3em  \hfil${##}$\hfil
      &\global\@eqcnt\tw@ \hskip0.35em
         $\displaystyle{##}$\hfil\tabskip\@centering
      &\global\@eqcnt\thr@@ \hb@xt@\z@\bgroup\hss##\egroup
         \tabskip\z@skip
      \cr
}
\makeatother

\def\C{\mathbbm{C}}
\def\R{\mathbbm{R}}

\def\Cn{\C^n}
\def\Cmn{\C^{m\times n}}
\def\Cmm{\C^{m\times m}}
\def\Cnn{\C^{n\times n}}
\def\Cnm{\C^{n\times m}}

\def\rats{{\cal R}}

\def\polyn{{\cal P}_{\kern-1pt n}}

\def\ORDER{{\cal O}}
\def\eps{\varepsilon}

\def\eop{{\rm e}}
\def\iop{{\rm i}}

\def\BA{{\bf A}}  \def\BAs{{\bf A}\kern-1pt}
\def\BB{{\bf B}} 
\def\BC{{\bf C}} 
 
\def\BE{{\bf E}}

\def\BH{{\bf H}} 
\def\BI{{\bf I}}

\def\BL{{\bf L}} 
\def\BM{{\bf M}}

\def\BP{{\bf P}}  \def\BPs{\BP\kern-1pt} 
 
\def\BR{{\bf R}} \def\Br{{\bf r}}

 \def\Bu{{\bf u}}
\def\BV{{\bf V}} \def\Bv{{\bf v}}
\def\BW{{\bf W}} \def\Bw{{\bf w}}
\def\BX{{\bf X}} \def\Bx{{\bf x}}
\def\BY{{\bf Y}} \def\By{{\bf y}}
 
\def\BLambda{\mbox{\boldmath$\Lambda$}}

\def\BSigma{\mbox{\boldmath$\Sigma$}}
\def\Bell{\boldsymbol{\ell}}
\def\BDelta{\boldsymbol{\Delta}}
\def\BGamma{\boldsymbol{\Gamma}}
\def\BTheta{\boldsymbol{\Theta}}
\def\BUp{\boldsymbol{\Upsilon}}

\def\CO{{\cal O}}

\def\LL{\mathbbm{L}}
\def\LL{{\rm I\kern-1.25pt L}}
\def\sLL{\LL_s}

\def\rank{{\rm rank}}

\def\LLij#1#2{{\Bv_{#1}^*\Br_{#2}^{} - \Bell_{#1}^* \Bw_{#2}^{} \over \mu_{#1}-\lambda_{#2}}}
\def\sLLij#1#2{{\mu_{#1}^{}\Bv_{#1}^*\Br_{#2}^{} - \lambda_{#2}^{}\Bell_{#1}^* \Bw_{#2}^{} \over \mu_{#1}-\lambda_{#2}}}

\def\zhat{\widehat{z}}

\def\wh#1{\widehat{#1}}
\def\Ahat{\wh{\BA}}
\def\Bhat{\wh{\BB}}
\def\Chat{\wh{\BC}}
\def\Ehat{\wh{\BE}}
\def\Hhat{\wh{\BH}}

\def\sv{s}   

\def\rho{\varrho}

\def\spec{\sigma}
\def\psa{\sigma_\eps}
\def\psagd{\sigma_\eps^{(\gamma,\delta)}}
\def\psaGD#1#2{\sigma_\eps^{(#1,#2)}}

\def\determinant{{\rm det}}
\def\Cauchy{{\boldsymbol {\mathcal C}}}

\makeindex             

\begin{document}

\title*{Pseudospectra of Loewner Matrix Pencils}
\author{Mark Embree and A. Cosmin Ioni\cb{t}\u{a}}
\institute{Mark Embree \at Department of Mathematics 
and Computational Modeling and Data Analytics Division, Academy of Integrated Science,
Virginia Tech, Blacksburg, VA 24061 USA, \email{embree@vt.edu}
\and A. Cosmin Ioni\cb{t}\u{a} \at The MathWorks Inc., Natick, MA 01760, USA,
\email{cionita@mathworks.com}}
%
%
\maketitle

\vspace*{0em}
\centerline{\emph{Dedicated to Thanos Antoulas.}}
\vspace*{2em}


\abstract{Loewner matrix pencils play a central role in the system realization theory of
Mayo and Antoulas, an important development in data-driven modeling.
The eigenvalues of these pencils reveal system poles.
How robust are the poles recovered via Loewner realization?
With several simple examples, we show how pseudospectra of Loewner pencils 
can be used to investigate the influence of interpolation point location and partitioning
on pole stability, the transient behavior of the realized system, 
and the effect of noisy measurement data.  
We include an algorithm to efficiently compute such pseudospectra
by exploiting Loewner structure.
}



\section{Introduction}

The landmark systems realization theory of Mayo and Antoulas~\cite{MA07}
constructs a dynamical system that interpolates tangential frequency domain measurements 
of a multi-input, multi-output system.
Central to this development is the \emph{matrix pencil}  $z\LL-\sLL$ 
composed of Loewner and shifted Loewner matrices $\LL$ and $\sLL$ 
that encode the interpolation data.
When this technique is used for exact system recovery 
(as opposed to data-driven model reduction),
the eigenvalues of this pencil match the poles of the
original dynamical system.
However other spectral properties, including the 
sensitivity of the eigenvalues to perturbation, can differ significantly;
indeed, they depend upon the location of the interpolation points in the 
complex plane relative to the system poles,
and how the poles are partitioned.
Since in many scenarios one uses $z\LL-\sLL$ to learn about
the original system, such subtle differences matter.

Pseudospectra are sets in the complex plane that contain the eigenvalues 
but provide additional insight about the sensitivity of those eigenvalues
to perturbation and the transient behavior of the underlying dynamical system.
While most often used to analyze single matrices, pseudospectral concepts
have been extended  to matrix pencils (generalized eigenvalue problems).  

This introductory note shows several ways to use pseudospectra 
to investigate spectral questions involving Loewner pencils
derived from system realization problems.  
Using simple examples, we explore the following questions.
\begin{itemize}
\item How do the locations of the interpolation points and their 
      partition into ``left'' and ``right'' points affect 
      the sensitivity of the eigenvalues of $z\LL-\sLL$?
\item Do solutions to the dynamical system $\LL \dot\Bx(t) = \sLL\Bx(t)$ 
      mimic solutions to the original system $\dot\Bx(t) = \BA\Bx(t)$, 
      especially in the transient regime? 
      Does this agreement depend on the interpolation points?
\item How do noisy measurements affect the eigenvalues of $z\LL-\sLL$?
\end{itemize}
We include an algorithm for computing pseudospectra of an 
$n$-dimensional  Loewner pencil in $\ORDER(n^2)$ operations, 
improving the $\ORDER(n^3)$ cost for generic matrix pencils;
the appendix gives a MATLAB implementation.

Throughout this note, we use $\sigma(\cdot)$ to denote the spectrum (eigenvalues)
of a matrix or matrix pencil, and $\|\cdot \|$ to denote the vector 2-norm
and the matrix norm it induces.  (All definitions here can readily
be adapted to other norms, as needed.  The algorithm, however, is 
designed for use with the 2-norm.)

\vspace*{-1.5em}

\section{Loewner realization theory in a nutshell}

\vspace*{-5pt}
We briefly summarize Loewner realization theory,
as developed by Mayo and Antoulas~\cite{MA07}; see also~\cite{ALI17}.
Consider the linear, time-invariant dynamical system
\begin{eqnarray*}
 \BE \dot\Bx(t) &=& \BA\Bx(t) + \BB \Bu(t) \\[.25em]
     \By(t) &=& \BC\Bx(t) 
\end{eqnarray*}
for $\BA, \BE \in\Cnn$, $\BB\in\C^{n\times m}$, and $\BC\in\C^{p\times n}$,
with which we associate, via the Laplace transform, the transfer function
$\BH(z) = \BC(z\BE-\BA)^{-1}\BB$.

Given \emph{tangential measurements} of $\BH(z)$ we seek to build a
realization of the system that interpolates the given data.
More precisely, consider 
the \emph{right interpolation data}\\[.5em]
\hspace*{2em} $\bullet$ distinct interpolation points $\lambda_1,\ldots, \lambda_\rho \in \C$;\\[.25em]
\hspace*{2em} $\bullet$ interpolation directions $\Br_1, \ldots, \Br_\rho \in \C^m$;\\[.25em]
\hspace*{2em} $\bullet$ function values $\Bw_1, \ldots, \Bw_\rho \in \C^p$;\\[.5em]
\noindent
and \emph{left interpolation data}\\[.5em]
\hspace*{2em} $\bullet$ distinct interpolation points $\mu_1,\ldots, \mu_\nu \in \C$; \\[.25em]
\hspace*{2em} $\bullet$ interpolation directions $\Bell_1, \ldots, \Bell_\nu \in \C^p$;\\[.25em]
\hspace*{2em} $\bullet$ function values $\Bv_1, \ldots, \Bv_\nu \in \C^m$.\\[.25em]
Assume the left and right interpolation points are disjoint 
$\{\lambda_i\}_{i=1}^\rho \cap \{\mu_j\}_{j=1}^\nu = \emptyset$;
indeed in our examples we will consider all the left and right points to be distinct.

The interpolation problem seeks matrices $\wh{\BA}, \wh{\BE}, \wh{\BB}, \wh{\BC}$ 
for which the transfer function $\wh{\BH}(z) = \wh{\BC}(z\wh{\BE}-\wh{\BA})^{-1}\wh{\BB}$ 
interpolates the data: for $i=1,\ldots, \rho$ and $j=1,\ldots, \nu$,
\[ \wh{\BH}(\lambda_i) \Br_i = \Bw_i, \qquad
   \Bell_j^*\wh{\BH}(\mu_j)  = \Bv_j^*.
\] 
Two structured matrices play a crucial role in the development of Mayo and Antoulas~\cite{MA07}.
From the data, construct the \emph{Loewner} 
and \emph{shifted Loewner} matrices 
\begin{equation}
\LL = \left[\begin{array}{ccc}
\LLij{1}{1} & \cdots & \LLij{1}{\rho} \\
 \vdots & \ddots & \vdots \\
\LLij{\nu}{1} & \cdots & \LLij{\nu}{\rho} \\
\end{array}\right], 
\quad
\sLL = \left[\begin{array}{ccc}
\sLLij{1}{1} & \cdots & \sLLij{1}{\rho} \\
 \vdots & \ddots & \vdots \\
\sLLij{\nu}{1} & \cdots & \sLLij{\nu}{\rho} \\
\end{array}\right]
\label{eq:LLsLL}
\end{equation}
i.e., the $(i,j)$ entries of these $\nu\times \rho$ matrices have the form
\[ (\LL)_{i,j} = {\Bv_i^*\Br_j^{} - \Bell_i^* \Bw_j^{} \over \mu_i-\lambda_j},
\qquad
({\sLL})_{i,j} = {\mu_i^{} \Bv_i^*\Br_j^{} - \lambda_j^{} \Bell_i^* \Bw_j^{} \over \mu_i-\lambda_j}.\]
Now collect the data into matrices.
The right interpolation points, directions, and data are stored in
\[ \BLambda 
   = \left[\begin{array}{ccc} \lambda_1 \\ & \ddots \\ & & \lambda_\rho
           \end{array}\right]
     \!\in\! \C^{\rho\times \rho}, 
   \quad
   \BR 
   = \left[\begin{array}{ccc} \\ \Br_1 & \cdots \Br_\rho \\  &
           \end{array}\right]
     \!\in\! \C^{m\times \rho},
    \quad
   \BW 
   = \left[\begin{array}{ccc} \\ \Bw_1 & \cdots \Bw_\rho \\ &
           \end{array}\right]
   \!\in\! \C^{p\times \rho},
\]
while the left interpolation points, directions, and data are stored in
\[ \BM 
   = \left[\begin{array}{ccc} \mu_1\\ & \ddots \\ & & \mu_\nu
           \end{array}\right]
    \!\in\! \C^{\nu\times \nu}, 
   \quad
   \BL 
   = \left[\begin{array}{ccc} \\ \Bell_1 & \cdots \Bell_\nu \\  &
           \end{array}\right]
   \!\in\! \C^{p\times \nu},
    \quad
   \BV 
   = \left[\begin{array}{ccc} \\ \Bv_1 & \cdots \Bv_\nu \\ &
           \end{array}\right]
   \!\in\! \C^{m\times \nu}.
\]

\subsection{Selecting and arranging interpolation points} \label{sec:syl}
As Mayo and Antoulas observe, Sylvester equations connect these matrices:
\begin{equation}
   \LL \BLambda - \BM \LL = \BL^*\BW-\BV^*\BR, \qquad
   \sLL\BLambda - \BM\sLL = \BL^*\BW\BLambda-\BM\BV^*\BR.
\label{eq:SylvesterLL}   
\end{equation}
Just using the dimensions of the components, note that 
\[
{\rm rank}(\BL^*\BW),\ 
{\rm rank}(\BV^*\BR),\ 
{\rm rank}(\BL^*\BW\BLambda),\ 
{\rm rank}(\BV^*\BR\BM)
\le 
\min\{\nu, \rho, m, p\}.\]
Thus for modest $m$ and $p$, the Sylvester equations~\eqref{eq:SylvesterLL}
must have low-rank right-hand sides.%
\footnote{For single-input, single-output (SISO) systems, $m=p=1$, so the rank of the
right-hand sides of the Sylvester equations cannot exceed two.
The same will apply for multi-input, multi-output systems with identical
left and right interpolation directions: $\Bell_i\equiv\Bell$ for all $i=1,\ldots,\nu$ and
$\Br_j\equiv\Br$ for all $j=1,\ldots,\rho$.} 
This situation often implies the rapid decay of singular values of
solutions to the Sylvester equation~\cite{ASZ02,BT17,Pen00a,Pen00b,Sab06}.
While this phenomenon is convenient in the context of balanced truncation
model reduction (enabling low-rank approximations to the controllability 
and observability Gramians), it is less welcome in the Loewner realization setting, 
where the rank of $\LL$ should reveal the order of the original system:
fast decay of the singular values of $\LL$ makes this rank ambiguous.
Since $\BLambda$ and $\BM$ are diagonal, they are normal matrices, and hence 
Theorem~2.1 of~\cite{BT17} gives
\begin{equation}\label{eq:syl}
 \frac{s_{qk+1}(\LL)}{s_1(\LL)}
 \le \inf_{\phi \in \rats_{k,k}} \frac{\max \{ |\phi(\lambda)|: \lambda \in \{\lambda_1, \ldots, \lambda_\rho\}\}}{\min \{ |\phi(\mu)|: \mu \in \{\mu_1, \ldots, \mu_\nu\}\}},
\end{equation}
where $s_j(\cdot)$ denotes the $j$th largest singular value,
$q = \rank(\BL^*\BW-\BV^*\BR)$, and $\rats_{k,k}$ denotes the set of irreducible 
rational functions whose numerators and denominators are polynomials of degree $k$ or less.%
\footnote{The same bound holds for $s_{qk+1}(\sLL)/s_1(\sLL)$ with 
$q = \rank(\BL^*\BW\BLambda-\BM\BV^*\BR)$.}
The right hand side of~\eqref{eq:syl} will be \emph{small} when there exists some $\phi\in\rats_{k,k}$
for which all $|\phi(\lambda_i)|$ are small, while all $|\phi(\mu_j)|$ are large:
a good separation of $\{\lambda_i\}$ from $\{\mu_j\}$ is thus \emph{sufficient} to ensure
the rapid decay of the singular values of $\LL$ and $\sLL$.
(Beckermann and Townsend give an explicit bound for the 
singular values of Loewner matrices when the interpolation points fall in disjoint
real intervals~\cite[Cor.~4.2]{BT17}.)

In our setting, where we want the singular values of $\LL$ and $\sLL$ to reveal the system's order
(without artificial decay of singular values as an accident of the arrangement of interpolation points),
it is \emph{necessary} for one $|\phi(\lambda_i)|$ to be at least the same size as the
smallest value of $|\phi(\mu_j)|$ for \emph{all} $\phi\in\rats_{k,k}$.
Roughly speaking, we want the left and right interpolation points to be 
close together (even interleaved).
While this arrangement is \emph{necessary} for slow decay of the singular values, 
it does not alone prevent such decay,
as evident in Figure~\ref{fig:siso_rank}
(since~\eqref{eq:syl} is only an upper bound).

Another heuristic, based on the Cauchy-like structure of $\LL$ and $\sLL$,
also suggests the left and right interpolation points
should be close together. 
Namely, $\LL$ and $\sLL$ are a more general
form of the Cauchy matrix $(\Cauchy)_{i,j} = 1/(\mu_i-\lambda_j)$, 
whose determinant has the elegant formula (e.g.,~\cite[p.~38]{HJ13})
\begin{equation}\label{eq:det}
\determinant(\Cauchy) = \frac{\prod_{1\leq i<j \leq n}(\mu_j-\mu_i)(\lambda_i-\lambda_j)}{\prod_{1\leq i\le j \leq n}(\mu_i-\lambda_j)}.
\end{equation}
It is an open question if $\determinant(\LL)$ and $\determinant(\sLL)$ have similarly elegant formulas.
Nevertheless, $\determinant(\LL)$ and $\determinant(\sLL)$ do have the same
denominator as $\determinant(\Cauchy)$
(which can be checked by recursively subtracting the first row from all other rows when computing the determinant).
This observation suggests that to avoid artificially small determinants for $\LL$ and $\sLL$ 
(which, up to sign, are the products of the singular values)
it is \emph{necessary} for the denominator of~\eqref{eq:det} to be small, and, thus, for the
left and right interpolation points to be close together.

In practice, we often start with initial
interpolation points $x_1,\ldots,x_{2n}$ that we want to partition 
into left and right interpolation points
to form $\LL$ and $\sLL$. Our analysis of~\eqref{eq:det} suggests a simple way to 
arrange the interpolation points such that the denominator of
$\determinant(\LL)$ and $\determinant(\sLL)$ is small:
relabel the points to satisfy
\begin{equation}\label{eq:arrange}
x_k = \underset{k \leq q \leq 2n}{\arg\min}~|x_{k-1}-x_q|, \quad {\rm for}~~ k = 2,\ldots,2n,
\end{equation}
\[
\{\mu_i\} = \{x_1,x_3,\ldots,x_{2n-1}\}, \quad
\{\lambda_j\} = \{x_2,x_4,\ldots,x_{2n}\}.
\]
The greedy reordering in~\eqref{eq:arrange} ensures that
$|x_k-x_{k-1}|$ is small and allows us to simply interleave the left and right interpolation
points. Moreover,
when $x_1,\ldots,x_{2n}$ are located on a line, the reordering in~\eqref{eq:arrange} simplifies to
directly interleaving $\mu_i$ and $\lambda_j$ and, thus, it can be skipped.
This ordering need not be optimal, as we do not visit
all possible combinations of $\mu_i-\lambda_j$; 
it simply seeks a partition that yields a large determinant (which must also depend on
the interpolation \emph{data}).
We note its simplicity, effectiveness, and efficiency (requiring only $\CO(n^2)$ operations).

\subsection{Construction of interpolants}

Throughout we make the fundamental assumptions that
for all $\wh{z} \in \{\lambda_i\}_{i=1}^\rho \cup \{\mu_j\}_{j=1}^\nu$,
\[ \rank@(\wh{z}@@@\LL-\sLL)
 = \rank\left(\left[\begin{array}{c}\LL  \\ \sLL\end{array}\right]\right)
 = \rank\left(\left[\begin{array}{cc}\LL & \sLL\end{array}\right]\right)
 =: r,\]
and we presume the underlying dynamical system is controllable and observable.

When $r = \nu=\rho$ is the order of the system, Mayo and Antoulas show that
the transfer function $\Hhat(z) := \Chat(z\Ehat-\Ahat)^{-1}\Bhat$ defined by
\begin{equation} \label{eq:Loew}
 \Ehat = -\LL@, \quad \Ahat = -\sLL, \quad \Bhat = \BV^*, \quad \Chat = \BW
\end{equation}
interpolates the $\rho+\nu$ data values.

When $r < \max(\nu,\rho)$, fix some $\wh{z}\in \{\lambda_i\}_{i=1}^\rho \cup \{\mu_j\}_{j=1}^\nu$
and compute the (economy-sized) singular value decomposition
\[ \wh{z}@@\sLL-\LL = \BY\BSigma\BX^*,\]
with $\BY\in\C^{\nu\times r}$, $\BSigma\in\R^{r\times r}$, and $\BX\in\C^{\rho\times r}$.
Then with 
\begin{equation} \label{eq:redLoew}
\Ehat = -\BY^*\LL\BX@, \quad \Ahat = -\BY^*\sLL\BX@, \quad \Bhat = \BY^*\BV^*, \quad \Chat = \BW\BX,
\end{equation}
$\Hhat (z) := \Chat(z\Ehat-\Ahat)^{-1}\Bhat$ 
defines an $r$th order system that interpolates the data.

\section{Pseudospectra for matrix pencils}

Though introduced decades earlier, in the 1990s pseudospectra 
emerged as a popular tool for analyzing 
the behavior of dynamical systems
(see, e.g.,~\cite{TTRD93}), 
eigenvalue perturbations
(see, e.g.,~\cite{CF96}),
and stability of uncertain 
linear time-invariant (LTI) systems
(see, e.g.,~\cite{HP90}).  

\begin{definition} \label{def:psa}
For a matrix $\BA\in\Cnn$ and $\eps>0$, the \emph{$\eps$-pseudospectrum} of $\BA$ is
\begin{equation}
\hspace*{-8pt}\psa(\BA) = \{\mbox{$z \in \C$ is an eigenvalue of $\BA+\BGamma$ for some $\BGamma\in\Cnn$ 
with $\|\BGamma\|< \eps$}\}.\label{eq:psadef}
\end{equation}
\end{definition}
For all $\eps>0$, $\psa(\BA)$ is a bounded, open subset 
of the complex plane that contains the eigenvalues of $\BA$.  
Definition~\ref{def:psa}  motivates pseudospectra via
eigenvalues of perturbed matrices.  A numerical analyst studying
accuracy of a backward stable eigenvalue algorithm might be concerned
with $\eps$ on the order of $n \|\BA\| \eps_{\rm mach}$, where
$\eps_{\rm mach}$ denotes the \emph{machine epsilon} 
for the floating point system~\cite{Ove01}.
An engineer or scientist might consider $\psa(\BA)$ for much larger $\eps$ 
values, corresponding to uncertainty in parameters or data that contribute
to the entries of $\BA$.

Via the singular value decomposition, one can show that~(\ref{eq:psadef})
is equivalent to
\begin{equation}
 \psa(\BA) = \{z\in\C: \|(z\BI-\BA)^{-1}\| > 1/\eps\};  \label{eq:psa2}
\end{equation}
see, e.g., \cite[chap.~2]{TE05}.
The presence of the resolvent $(z\BI-\BA)^{-1}$ in this definition
suggests a connection to the transfer function $\BH(z) = \BC(z\BI-\BA)^{-1}\BB$
for the system
\begin{eqnarray*}
\dot\Bx(t) &=& \BA\Bx(t) + \BB \Bu(t) \\[.25em]
\By(t) &=& \BC\Bx(t).
\end{eqnarray*}
Indeed, definition~\eqref{def:psa} readily leads to bounds on $\|\eop^{t\BA}\|$,
and hence transient growth of solutions to $\dot\Bx(t) = \BA\Bx(t)$; 
see \cite[part~IV]{TE05}.

Various extensions of 
pseudospectra have been proposed
to handle more general eigenvalue problems and dynamical systems;
see~\cite{Emb14} for a concise survey.
The first elaborations addressed the generalized eigenvalue problem
$\BA\Bx = \lambda \BE\Bx$ 
(i.e., the matrix pencil $z@\BE-\BA$)~\cite{FGNT96,Rie94,Ruh95}.
Here we focus on the definition proposed by
Frayss\'e, Gueury, Nicoud, and Toumazou~\cite{FGNT96}, 
which is ideally suited to analyzing eigenvalues of nearby matrix 
pencils.
To permit the perturbations to $\BA$ and $\BE$ to be scaled independently,
this definition includes two additional parameters, $\gamma$ and $\delta$. 

\begin{definition} \label{def:psagd}
Let $\gamma,\delta> 0$.  For a pair of matrices $\BA,\BE\in\Cnn$
and any $\eps>0$, 
the \emph{$\eps$-$(\gamma,\delta)$-pseudospectrum} $\psagd(\BA,\BE)$ 
of the matrix pencil $z\BE - \BA$ is the set
\begin{eqnarray*} \label{eq:psaAE}
\psagd(\BA,\BE) &=& \{\mbox{$z \in \C$ is an eigenvalue of the pencil 
$z(\BE+\BDelta) - (\BA+\BGamma)$} \label{eq:psa2} \\
&& \ \ \ \ \ \mbox{for some $\BGamma, \BDelta \in \Cnn$ with 
$\|\BGamma\|< \eps \gamma$, $\|\BDelta\|< \eps\delta$}\}.  \nonumber
\end{eqnarray*}
\end{definition}
This definition has been extended to matrix polynomials in~\cite{HT02,TH01}.
\begin{remark} \label{rem:unbdd}
Note that $\psagd(\BA,\BE)$ is an open, nonempty subset of the complex plane,
but it need not be bounded.  
\begin{enumerate}
\item[(a)] If $z@\BE-\BA$ is a singular pencil ($\rank(z@\BE-\BA)<n$ for all $z\in\C$), then $\sigma(\BA,\BE) = \C$.
\item[(b)] If $z@\BE-\BA$ is a regular pencil but $\BE$ is not invertible, 
then $\sigma(\BA,\BE)$ contains an infinite eigenvalue.
\item[(c)] If $\BE$ is nonsingular but $\eps \delta$ exceeds the distance of $\BE$ 
to singularity (the smallest singular value of $\BE$), then
$\psagd(\BA,\BE)$ contains the point at infinity.
\end{enumerate}
Since these pseudospectra can be unbounded,
Lavall\'ee~\cite{Lav97} and Higham and Tisseur~\cite{HT02} visualize $\psagd(\BA,\BE)$
as stereographic projections on the Riemann sphere.
\end{remark}

\begin{remark}
Just as the conventional pseudospectrum $\psa(\BA)$ can be characterized
using the resolvent of $\BA$ in~(\ref{eq:psa2}), 
Frayse\'{e} et al.~\cite{FGNT96} show that Definition~\ref{def:psagd}
is equivalent to
\begin{equation}
 \psagd(\BA,\BE) = \left\{ z\in\C: \|(z\BE-\BA)^{-1}\| > {1\over \eps(\gamma + |z| \delta)}\right\}.
\label{eq:psagd2}
\end{equation}
This formula suggests a way to compute $\psagd(\BA,\BE)$:  
evaluate $ \|(z\BE-\BA)^{-1}\|$ on a grid of points covering a relevant
region of the complex plane (or the Riemann sphere) and use a contour plotting
routine to draw boundaries of $\psagd(\BA,\BE)$.\ \   
The accuracy of the resulting pseudospectra depends on the density of the grid.
Expedient algorithms for computing $\|(z\BE-\BA)^{-1}\|$ can be derived by
computing a unitary simultaneous triangularization (generalized Schur form) 
of $\BA$ and $\BE$ in $\ORDER(n^3)$ operations, then using inverse iteration 
or inverse Lanczos, as described by Trefethen~\cite{Tre99a} and 
Wright~\cite{Wri02b}, to compute $\|(z\BE-\BA)^{-1}\|$
at each grid point in $\ORDER(n^2)$ operations.  
For the structured Loewner pencils of interest here, one can compute
 $\|(z\BE-\BA)^{-1}\|$ in $\ORDER(n^2)$ operations without recourse to the 
$\ORDER(n^3)$ preprocessing step, as proposed in Section~\ref{sec:fastiter}.
\end{remark}

\begin{remark}
Definition~\ref{def:psagd} can be extended to $\delta=0$ by 
only perturbing $\BA$:
\begin{eqnarray*} \label{eq:psaAE0}
\psaGD{1}{0}(\BA,\BE) &=& \{\mbox{$z \in \C$ is an eigenvalue of the pencil 
$z@\BE-(\BA+\BGamma)$} \label{eq:psa2} \\
&& \ \ \ \ \ \mbox{for some $\BGamma\in\Cnn$ with $\|\BGamma\|< \eps$}\}.  \nonumber \\
&=& \{ z\in\C: \|(z\BE-\BA)^{-1}\| > 1/\eps \}. 
\end{eqnarray*}
This definition may be more suitable for cases where $\BE$ is fixed and 
uncertainty in the system only emerges, e.g., through physical parameters
that appear in $\BA$.
\end{remark}

\begin{remark} 
Since we ultimately intend to study the pseudospectra 
$\psagd(\sLL,\LL)$ of Loewner matrix pencils,
one might question the use of generic perturbations
$\BGamma, \BDelta\in\Cnn$ in Definition~\ref{def:psagd}.
Should we restrict $\BGamma$ and $\BDelta$ to maintain
Loewner structure, 
i.e., so that $\sLL+\BGamma$ and $\LL+\BDelta$ maintain
the coupled shifted Loewner--Loewner form 
in~(\ref{eq:LLsLL})?
Such sets are called \emph{structured pseudospectra}.

Three considerations motivate the study of generic perturbations
$\BGamma, \BDelta\in\Cnn$: one practical, one speculative, and one philosophical.
(a)~Beyond repeatedly computing the eigenvalues of Loewner pencils with
randomly perturbed data, no systematic way is known to compute
the coupled Loewner structured pseudospectra, 
i.e., no analogue of the resolvent-like definition~(\ref{eq:psagd2})
is known.
(b)~Rump~\cite{Rum06} showed that in many cases,
preserving structure has little effect on the standard 
matrix pseudospectra.  For example, the structured 
$\eps$-pseudospectrum of a Hankel matrix $\BH$ allowing only
complex Hankel perturbations exactly matches the
unstructured $\eps$-pseudospectrum $\psa(\BH)$ 
based on generic complex perturbations~\cite[Thm.~4.3]{Rum06}.
\emph{Whether a similar results holds for Loewner structured pencils is an interesting open question.}
(c)~If one seeks to analyze the \emph{behavior} of
dynamical systems (as opposed to eigenvalues of nearby matrices),
then generic perturbations give much greater insight;
see~\cite[p.~456]{TE05} for an example where \emph{real-valued} 
perturbations do not move the eigenvalues much toward the imaginary axis
(hence the real structured pseudospectra are benign),
yet the stable system still exhibits strong transient growth.
\end{remark}

%


As we shall see in Section~\ref{sec:fullrank}, Definition~\ref{def:psagd} provides a 
helpful tool for investigating the sensitivity of eigenvalues of matrix pencils.
A different generalization of Definition~\ref{def:psa} 
gives insight into the transient behavior 
of solutions of $\BE\dot\Bx(t)=\BA\Bx(t)$.
This approach is discussed in~\cite[chap.~45]{TE05}, following~\cite{Rie94,Ruh95},
and has been extended to handle singular $\BE$ in~\cite{EK17} 
(for differential-algebraic equations and descriptor systems).
Restricting our attention here to nonsingular $\BE$,
we analyze the conventional (single matrix) pseudospectra $\sigma_\eps(\BE^{-1}\BA)$.  
From these sets one can develop various upper and lower bounds on 
$\|\eop^{t\BE^{-1}\BA}\|$ and $\|\Bx(t)\|$~\cite[chap.~15]{TE05}.
Here we shall just state one basic result.
If $\sup \{{\rm Re}(z) : z\in\sigma_\eps(\BE^{-1}\BA)\} = K \eps$ for some $K\ge 1$
(where ${\rm Re}(\cdot)$ denotes the real part of a complex number),
then
\begin{equation} \label{eq:kreiss}
\sup_{t\ge 0} \|\eop^{t@\BE^{-1}\BA}\| \ge K.
\end{equation}
This statement implies that there exists some unit-length initial condition
$\Bx(0)$ such that $\|\Bx(t)\| \ge K$, even though $\spec(\BA,\BE)$ may be contained
in the left half-plane.
(Optimizing this bound over $\eps>0$ yields the Kreiss Matrix Theorem~\cite[(15.9)]{TE05}.)

\medskip
Pseudospectra of matrix pencils provide a natural vehicle to explore that stability of the
matrix pencil associated with the Loewner realization in~(\ref{eq:Loew}).  
We shall thus investigate eigenvalue perturbations via $\psagd(\Ahat,\Ehat) = \psagd(\sLL,\LL)$ 
and transient behavior via $\psa(\LL^{-1}\sLL)$.

\section{Efficient computation of Loewner pseudospectra} \label{sec:fastiter}

We first present a novel technique for efficiently computing pseudospectra of 
large Loewner matrix pencils, $\psagd(\sLL,\LL)$,
using the equivalent definition given in~(\ref{eq:psagd2}).
When the Loewner matrix $\LL$ is nonsingular, we employ inverse iteration to
exploit the structure of the Loewner pencil to compute $\|(z\LL-\sLL)^{-1}\|$ (in the two-norm)
using only $\CO(n^2)$ operations.  This avoids the need to  compute an initial simultaneous unitary
triangularization of $\sLL$ and $\LL$ using the QZ algorithm, an $\CO(n^3)$ operation.

Inverse iteration (and inverse Lanczos) for $\|(z\LL-\sLL)^{-1}\|$ requires computing
\begin{equation}
(z\LL-\sLL)^{-*}(z\LL-\sLL)^{-1}\Bu
\label{eq:inviter}
\end{equation}
for a series of vectors $\Bu \in \Cn$ (e.g., see \cite[chap.~39]{TE05}).
We invoke a property observed by Mayo and Antoulas~\cite{MA07}, related to~\eqref{eq:SylvesterLL}: 
by construction, the Loewner and
shifted Loewner matrices satisfy
\[
\sLL - \LL\BLambda = \BV^*\BR, \qquad \sLL - \BM\LL = \BL^*\BW.
\]
Thus the resolvent can be expressed using only $\LL$ and not $\sLL$:
\[
 (z\LL-\sLL)^{-1} = (\LL(z\BI-\BLambda)-\BV^*\BR)^{-1}.
\]
We now use the Sherman--Morrison--Woodbury formula (see, e.g.,~\cite{GV12}) to get
\[
 (z\LL-\sLL)^{-1}
 = \left( \BI + \BUp(z) (\BI - \BR\BUp(z))^{-1} \BR \right) (z\BI-\BLambda)^{-1} \LL^{-1}
 = \BTheta(z)\LL^{-1},
\]
where $\BUp(z) := (z\BI-\BLambda)^{-1} \LL^{-1}\BV^*$ 
and $\BTheta(z) := \left(\BI + \BUp(z) (\BI - \BR\BUp(z))^{-1} \BR \right) (z\BI-\BLambda)^{-1}$. 
As a result, we can compute the inverse iteration
vectors in (\ref{eq:inviter}) as
\begin{equation}
(z\LL-\sLL)^{-*}(z\LL-\sLL)^{-1}\Bu = \LL^{-*}\BTheta(z)^*\BTheta(z)\LL^{-1}\Bu,
\label{eq:inviter2}
\end{equation}
which requires solving several linear systems given by the same
Loewner matrix $\LL$, e.g., $\LL^{-1}\Bu$, $\LL^{-1}\BV^*$.

Crucially, solving a linear system involving a Loewner matrix $\LL \in \Cnn$ can be done efficiently in
only $2(m+p+1)n^2$ operations because $\LL$ has displacement rank $m+p$.
More precisely, $\LL$ is a Cauchy-like matrix that
satisfies the Sylvester equation~(\ref{eq:SylvesterLL})
given by diagonal generator matrices $\BLambda$ and $\BM$ and a right-hand side
of rank at most $m+p$, i.e., ${\rm rank}(\BL^*\BW-\BV^*\BR) \le m+p$. The displacement rank structure of
$\LL$ can be exploited to compute its LU factorization in only $2(m+p)n^2$
flops (see \cite{GKO95b} and \cite[sect.~12.1]{GV12}). Given the LU factorization of $\LL$,
solving $\LL^{-1}\Bu$ in (\ref{eq:inviter2}) via standard forward and backward substitution
requires another $2n^2$ operations.

Next, multiplying $\BTheta(z)$ with the solution of $\LL^{-1}\Bu$ requires a total of
$2mn^2 + (4m^2+6m+2)n + \frac{2}{3}m^3-m^2$ operations, namely (to leading order on the factorizations):\\[.5em]
\hspace*{2em} $\bullet$ $2mn^2 + 2mn$ operations to compute $\BUp(z) \in \Cnm$;\\[.25em]
\hspace*{2em} $\bullet$ $m^2(2n-1) + m$ operations to compute $\BI - \BR\BUp(z) \in \Cmm$;\\[.25em]
\hspace*{2em} $\bullet$ $\frac{2}{3}m^3+2m^2n$ operations to solve
$(\BI - \BR\BUp(z))^{-1} \BR \in \Cmn$ via \\
\hspace*{3em} an LU factorization followed by $n$ forward and backward substitutions;\\[.25em]
\hspace*{2em} $\bullet$ $n + m(2n-1) + (2m-1)n + 2n$ operations to multiply $\BTheta(z)$ with $\LL^{-1}\Bu$.\\[.5em]
\indent
Finally, multiplying with $\LL^{-*}\BTheta(z)^*$ in~(\ref{eq:inviter2}) requires an additional 
$2n^2 + 2mn^2 + (4m^2+6m+2)n + \frac{2}{3}m^3-m^2$ operations, bringing the total cost of 
computing (\ref{eq:inviter2}) to 
$2(3m+p+1)n^2 + 4(2m^2+3m+1)n + \frac{4}{3}m^3-2m^2$ operations.

In practice, the sizes of the right 
and left tangential directions are much smaller than the size of the Loewner pencil, i.e.,
$m, p \ll n$. For example, for scalar data (associated with SISO systems), $m = p = 1$. 
Therefore, in practice, computing
(\ref{eq:inviter2}) can be done in only $\CO(n^2)$ operations.

Partial pivoting can be included in the LU factorization of the
Loewner matrix $\LL$ to overcome numerical difficulties. Adding partial pivoting maintains the
$\CO(n^2)$ operation count for the LU factorization of $\LL$ (see \cite[sect.~12.1]{GV12}),
and hence computing (\ref{eq:inviter2}) can still be done in $\CO(n^2)$ operations.
The appendix gives a MATLAB implementation of this efficient inverse iteration.

We measure these performance gains for a Loewner pencil generated by sampling
$f(x) = \sum_{k=1}^8 (-1)^{k+1}\left(1+100(x-k)^2\right)^{-1/2}
+ (-1)^{k+1}\left(1+100(x-k-1/2)^2\right)^{-1/2}$
at $2n$ points uniformly spaced in the interval $[1,8]$. We compare our new $\CO(n^2)$ 
Loewner pencil inverse iteration against a standard implementation
(see \cite[p.~373]{TE05}) applied to a simultaneous triangularization of $\sLL$ and $\LL$. 
(The simultaneous triangularization costs $\CO(n^3)$ but is fast, as MATLAB's {\tt qz} routine 
invokes LAPACK code.  For a fair comparison, we test against a C++ implementation of the
fast Loewner code, compiled into a MATLAB {\tt .mex} file.)
Table~\ref{tab:fast} shows timings for both implementations by computing $\|(z\LL-\sLL)^{-1}\|$
on a $200\times 200$ grid of points.
As expected, exploiting the Loewner structure gives a significant performance improvement for large $n$.

\begin{table}[t!]
\caption{\label{tab:fast}
Comparison of speed of computing $\psa^{(1,1)}(\sLL,\LL)$ using the
fast Loewner algorithm versus a generic inverse iteration method applied to
simultaneous triangularizations of $\sLL$ and $\LL$.}
\begin{center}
\begin{tabular}{ccc|crc|crc}
&$n$&&&$\CO(n^2)$ algorithm &&&$\CO(n^3)$ algorithm &\\ \hline \hline
&$100$&&&$1.85$ seconds\ \ \ \ \ &&&$1.65$ seconds\ \ \ \ \ &\\\hline
&$200$&&&$6.75$ seconds\ \ \ \ \ &&&$6.75$ seconds\ \ \ \ \ &\\\hline
&$300$&&&$17.24$ seconds\ \ \ \ \ &&&$33.30$ seconds\ \ \ \ \ &\\\hline
&$400$&&&$49.15$ seconds\ \ \ \ \ &&&$65.05$ seconds\ \ \ \ \ &\\\hline
\end{tabular}
\end{center}
\end{table}

We next examine two simple examples involving full-rank realization of SISO systems,
to illustrate the kinds of insights one can draw from pseudospectra of Loewner pencils.


\section{Example 1: eigenvalue sensitivity and transient behavior}

We first consider a simple controllable and observable SISO system with $n=2$:
\begin{equation}
\BA = \left[\begin{array}{cc} -1.1 & 1 \\ 1 & -1.1 \end{array}\right], \qquad
   \BB = \left[\begin{array}{c} 0 \\ 1 \end{array}\right], \qquad 
   \BC = \left[\begin{array}{cc} 0 & 1 \end{array}\right].  \label{eq:twodim1}
\end{equation}
This $\BA$ is symmetric negative definite, with eigenvalues $\sigma(\BA) = \{-0.1,-2.1\}$.
Since the system is SISO, the transfer function $\BH(s) = \BC(s\BI-\BA)^{-1}\BB$ 
maps $\C$ to $\C$, and hence the choice of ``interpolation directions'' is trivial
(though the division into ``left'' and ``right'' points matters).
We take $\rho=\nu=2$ left and right interpolation points, 
with $\Br_1 = \Br_2 = 1$ and $\Bell_1 = \Bell_2 = 1$.
We will study various choices of interpolation points, all of which satisfy,
for each $\zhat \in \{\lambda_1, \lambda_2, \mu_1, \mu_2\}$,
\begin{equation} 
\rank(\zhat\, \LL - \sLL) = \rank(\LL) = \rank(\sLL) = n=2.\label{eq:Lrank2}
\end{equation}
This basic set-up makes it easy to focus on the influence of the interpolation points $\lambda_1$, $\lambda_2$, $\mu_1$, $\mu_2$. 
We will use the pseudospectra $\psa^{(1,1)}(\sLL,\LL)$ to examine how the interpolation points 
affect the stability of the eigenvalues of the Loewner pencil.

\begin{table}[b!]
\caption{\label{tbl:twodim1}Right and left interpolation points for the two-dimensional
SISO system~\eqref{eq:twodim1}.  The two right columns report the singular values
of the Loewner matrix $\LL$.}

\begin{center}
\begin{tabular}{c|crcrcrcrc|clcl}
\emph{example}& \ &
\multicolumn{1}{c}{$\lambda_1$} & \ \ & 
\multicolumn{1}{c}{$\lambda_2$} & \ \ & 
\multicolumn{1}{c}{$\mu_1$} & \ \ & 
\multicolumn{1}{c}{$\mu_2$} & \ \ &  \ \ &
\multicolumn{1}{c}{$s_1(\LL)$} & \ \ & 
\multicolumn{1}{c}{$s_2(\LL)$} \\ \hline \hline
(a) & & \rule{0pt}{9pt} 0 && 1 && $1@\iop$ && $-1@\iop$ &&&  6.9871212 && 0.0731542 \\ \hline
(b) & & 0.25 && 0.75 && $2@\iop$  &&  $-2@\iop$ &&& 1.0021659 && 0.0296996 \\ \hline
(c) & & 0.40 && 0.60 && $4@\iop$ && $-4@\iop$ &&& 0.3605151 && 0.0057490 \\ \hline
(d) & & 8 && 9 && 10 && 11 &&& 0.0035344 && 0.0000019 \\ \hline
\end{tabular}
\end{center}
\end{table}

\begin{figure}[b!]
\begin{center}
    \includegraphics[width=2.15in]{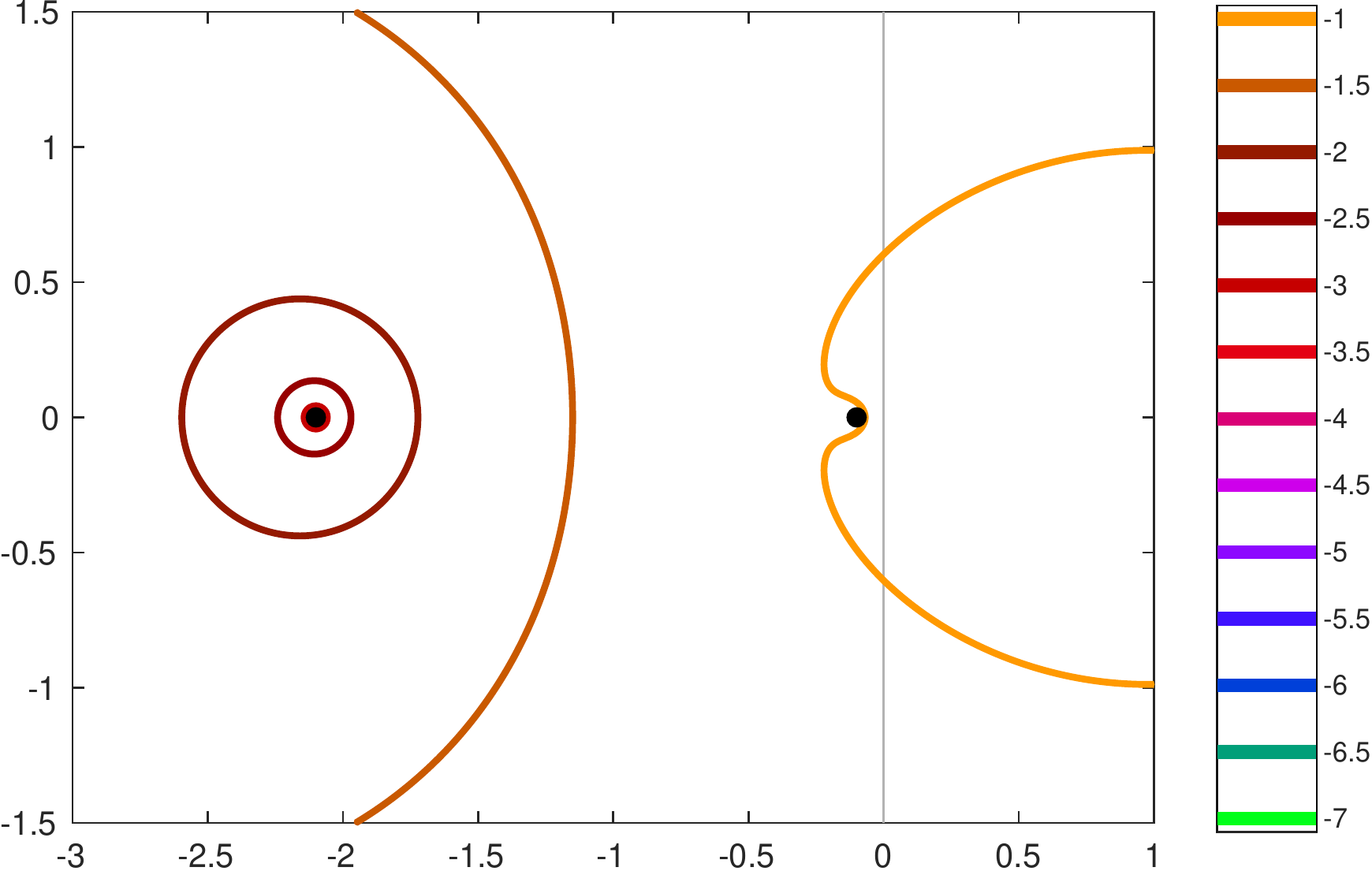} \quad     \includegraphics[width=2.15in]{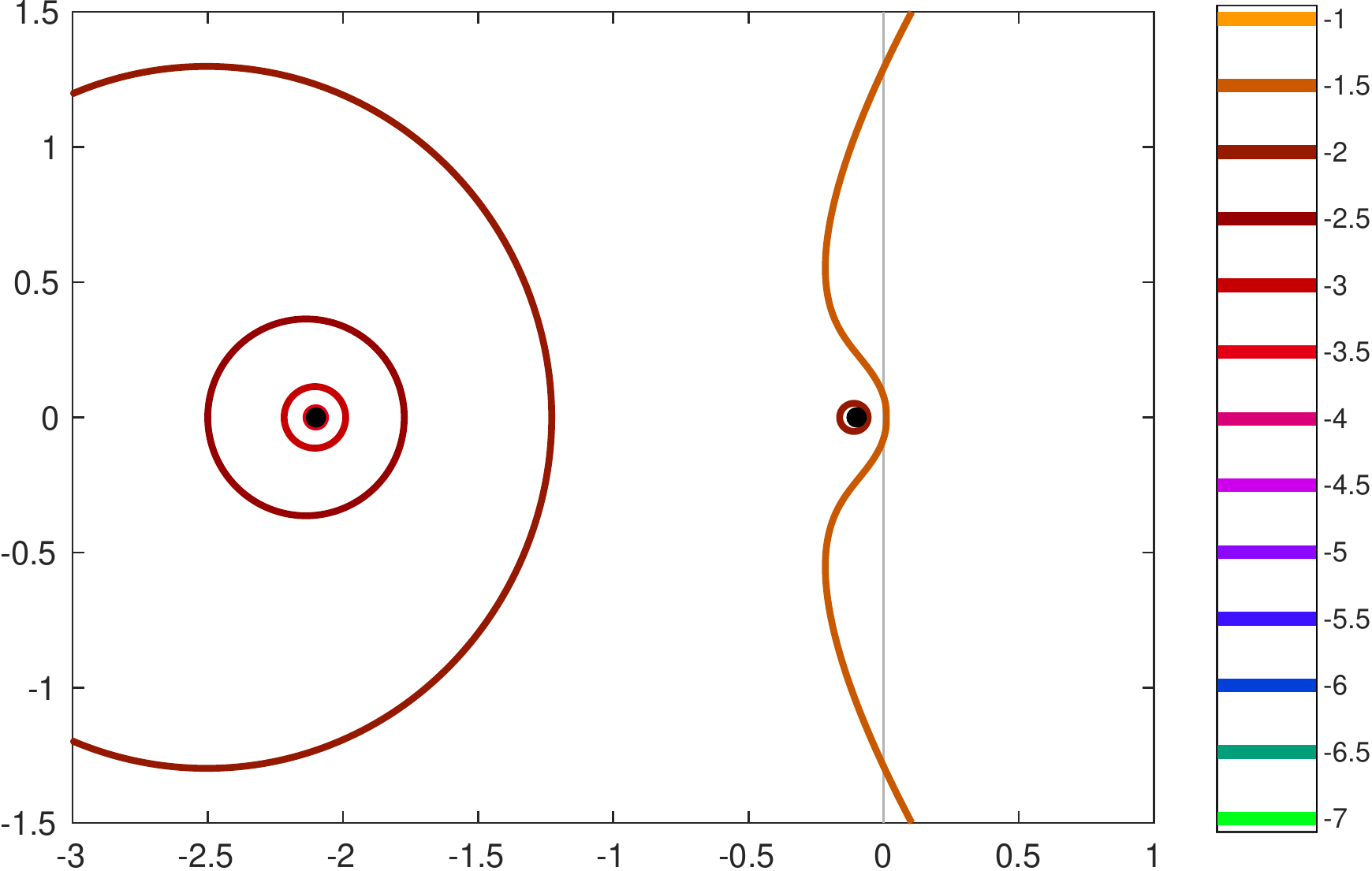}
    \begin{picture}(0,0)
       \put(-235,12){(a)}
       \put(-72,12){(b)}
    \end{picture}

    \vspace*{5pt}
    \includegraphics[width=2.15in]{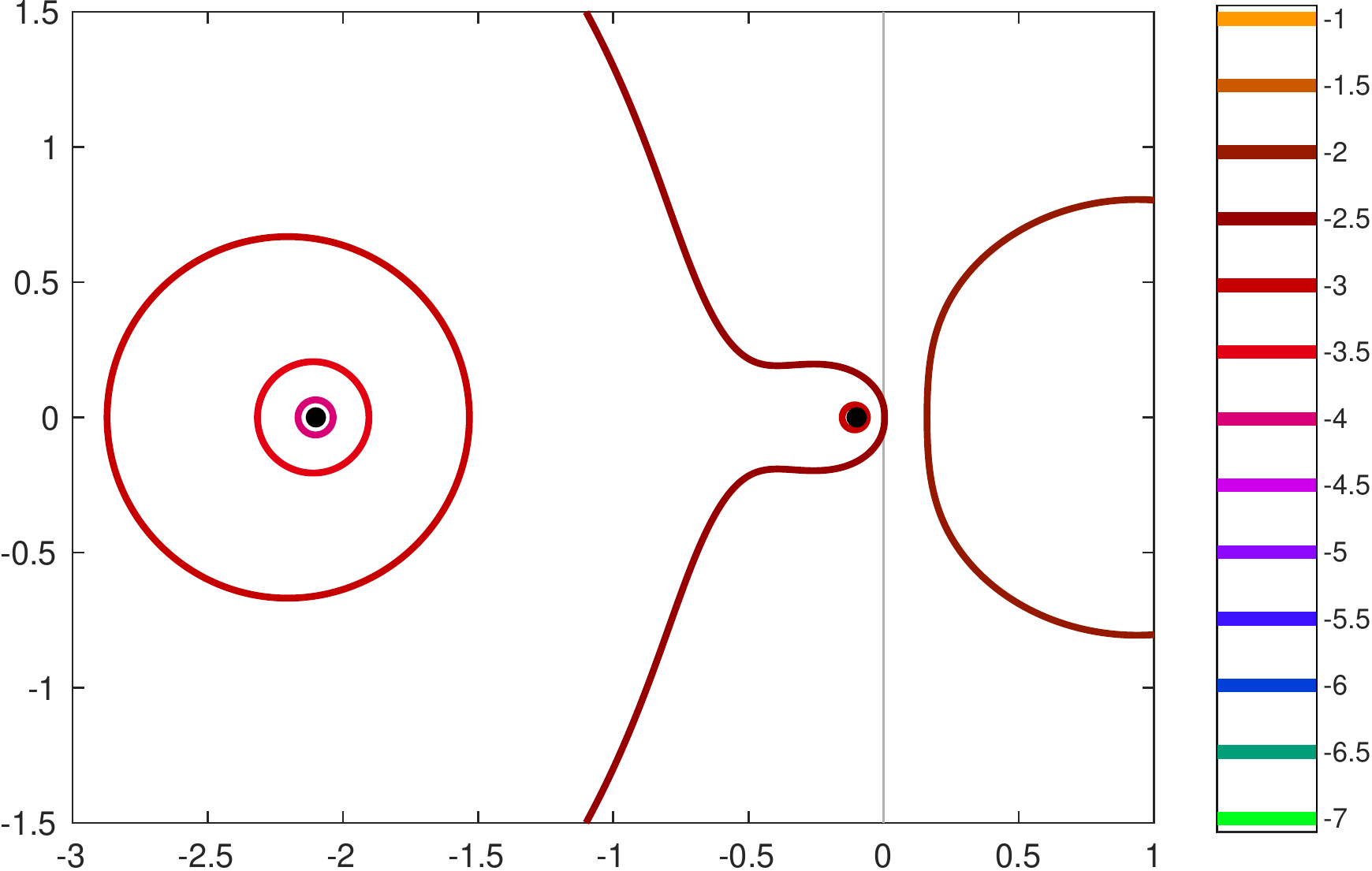} \quad     \includegraphics[width=2.15in]{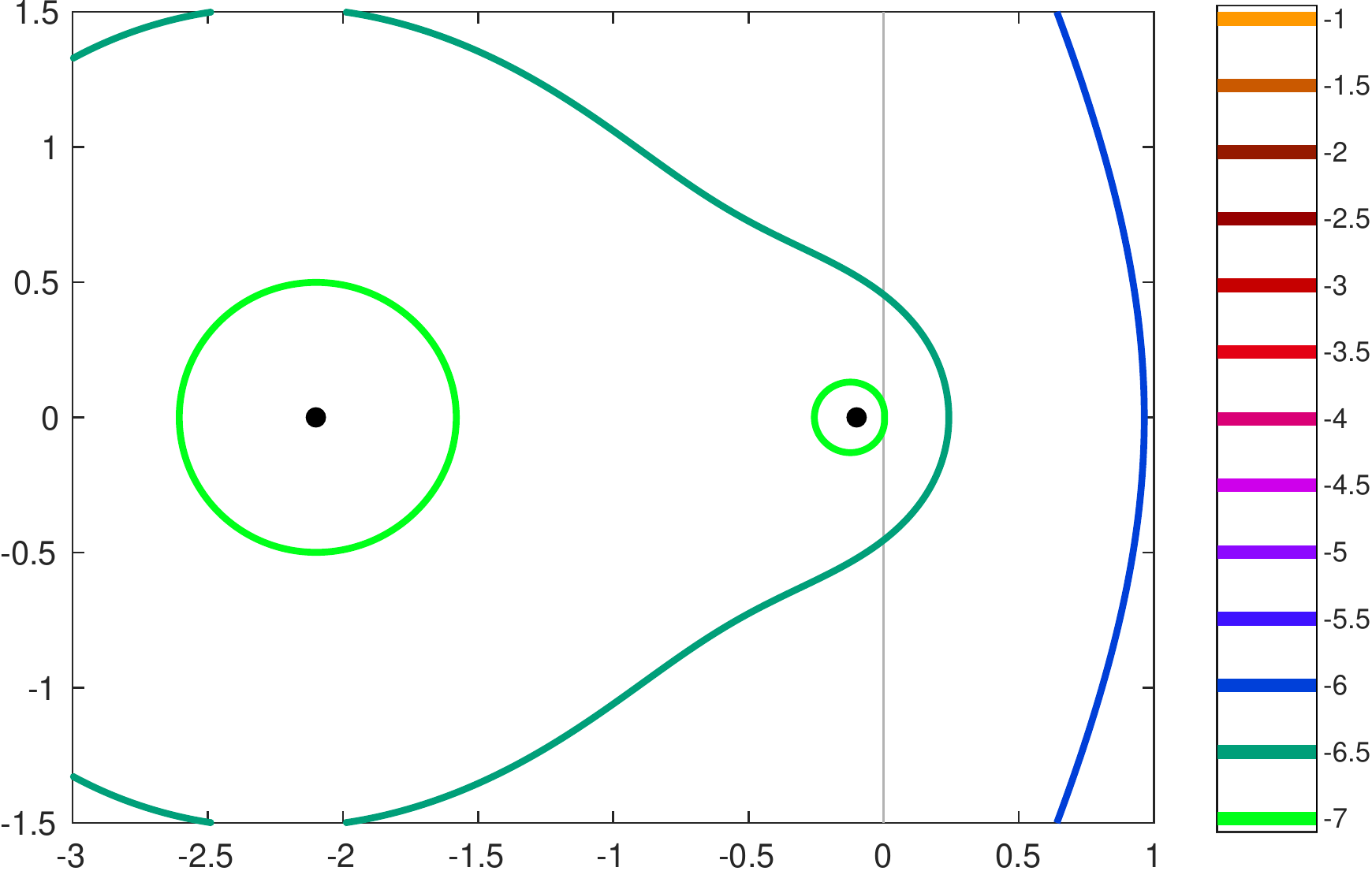}
    \begin{picture}(0,0)
       \put(-235,12){(c)}
       \put(-72,12){(d)}
    \end{picture}
\end{center}

\vspace*{-5pt}
    \caption{\label{fig:twodim1}
    Boundaries of pseudospectra $\psa^{(1,1)}(\sLL,\LL)$ for four Loewner realizations of the 
    system~\eqref{eq:twodim1} using the interpolation points in Table~\ref{tbl:twodim1}.
    All four realizations correctly give $\spec(\sLL,\LL) = \spec(\BA)$, but 
    $\psa^{(1,1)}(\sLL,\LL)$ show how the stability of the realized eigenvalues 
    depends on the choice of interpolation points. 
    \emph{In this and all similar plots, the colors denote $\log_{10}(\eps)$.}
    Thus, in plot~(d), there exist perturbations to $\sLL$ and $\LL$ of norm $10^{-6.5}$ 
    that move an eigenvalue into the right half-plane.
 }
\end{figure}

Table~\ref{tbl:twodim1} records four different choices of $\{\lambda_1, \lambda_2, \mu_1, \mu_2\}$;
Figure~\ref{fig:twodim1} shows the corresponding pseudospectra $\psa^{(1,1)}(\sLL,\LL)$.
All four Loewner realizations match the eigenvalues of $\BA$
and satisfy the interpolation conditions.  
However, the pseudospectra show how the \emph{stability} of the eigenvalues 
$-0.1$ and $-2.1$ differs across these four realizations.  
These eigenvalues become increasingly sensitive from example~(a) to~(d),
as the interpolation points move farther from $\sigma(\BA)$.
Table~\ref{tbl:twodim1} also shows the singular values of
the Loewner matrix $\LL$, demonstrating how the second singular value $s_2(\LL)$
decreases as the eigenvalues become increasingly sensitive.  
(Taken to a greater extreme, it would eventually be difficult to determine
if $\LL$ truly is rank~2.)

\begin{remark}
By Remark~\ref{rem:unbdd} on page~\pageref{rem:unbdd}, 
note that if $\eps< \sv_{\rm min}(\LL)$, then $\psaGD{1}{1}(\sLL,\LL)$ will be
unbounded.  Thus the decreasing values of $\sv_{\rm min}(\LL)$ in Table~\ref{tbl:twodim1}
suggests the enlarging pseudospectra seen in Figure~\ref{fig:twodim1}.
For example, in case~(d) the $\eps=10^{-5}$ pseudospectrum $\psaGD{1}{1}(\sLL,\LL)$
must contain the point at infinity.
\end{remark}

\begin{figure}[p]
\begin{center}
    \includegraphics[width=2.15in]{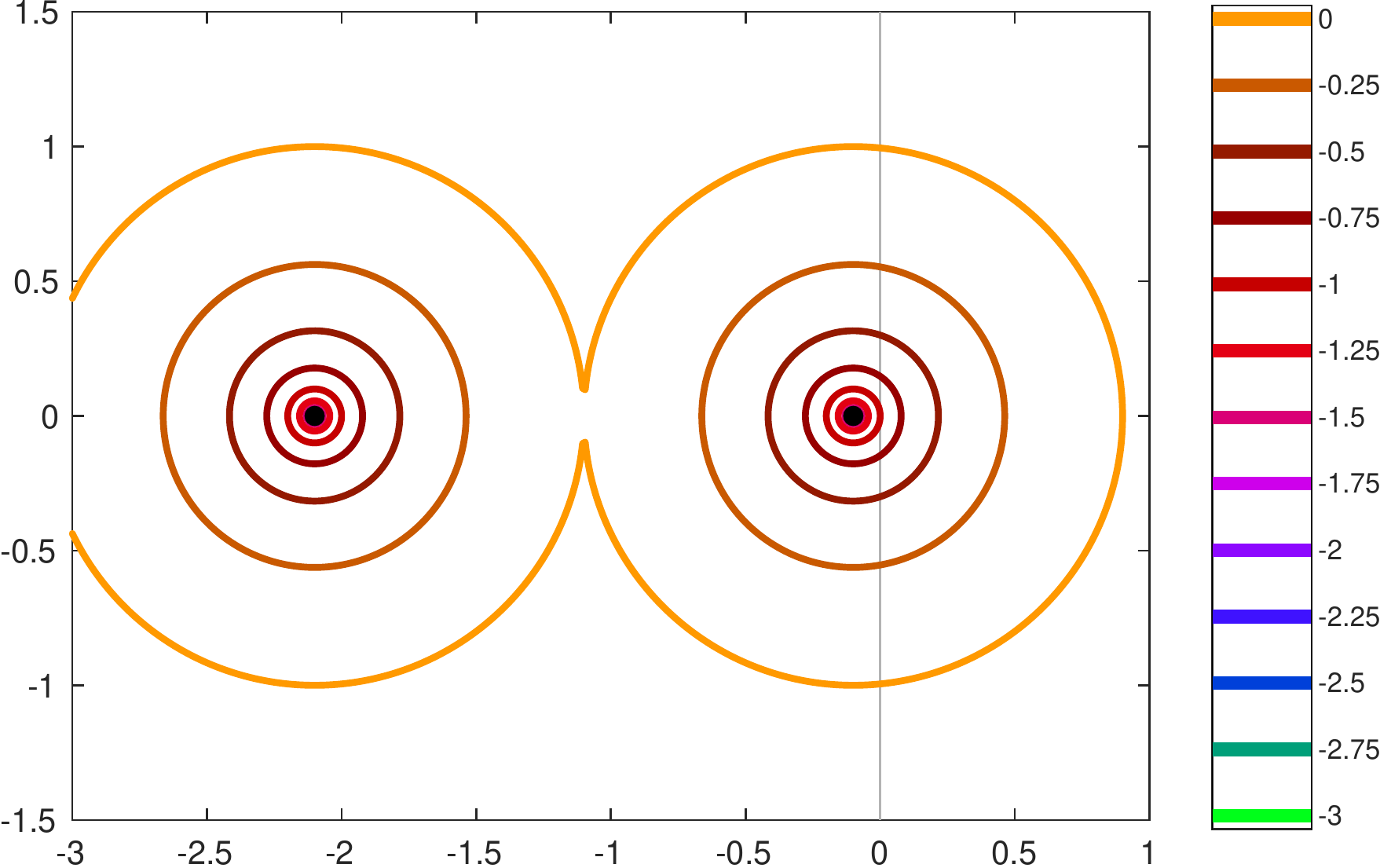}
    \begin{picture}(0,0)
       \put(-55,11){$\psa(\BA)$}
    \end{picture}
    
    \vspace*{7pt}

    \includegraphics[width=2.15in]{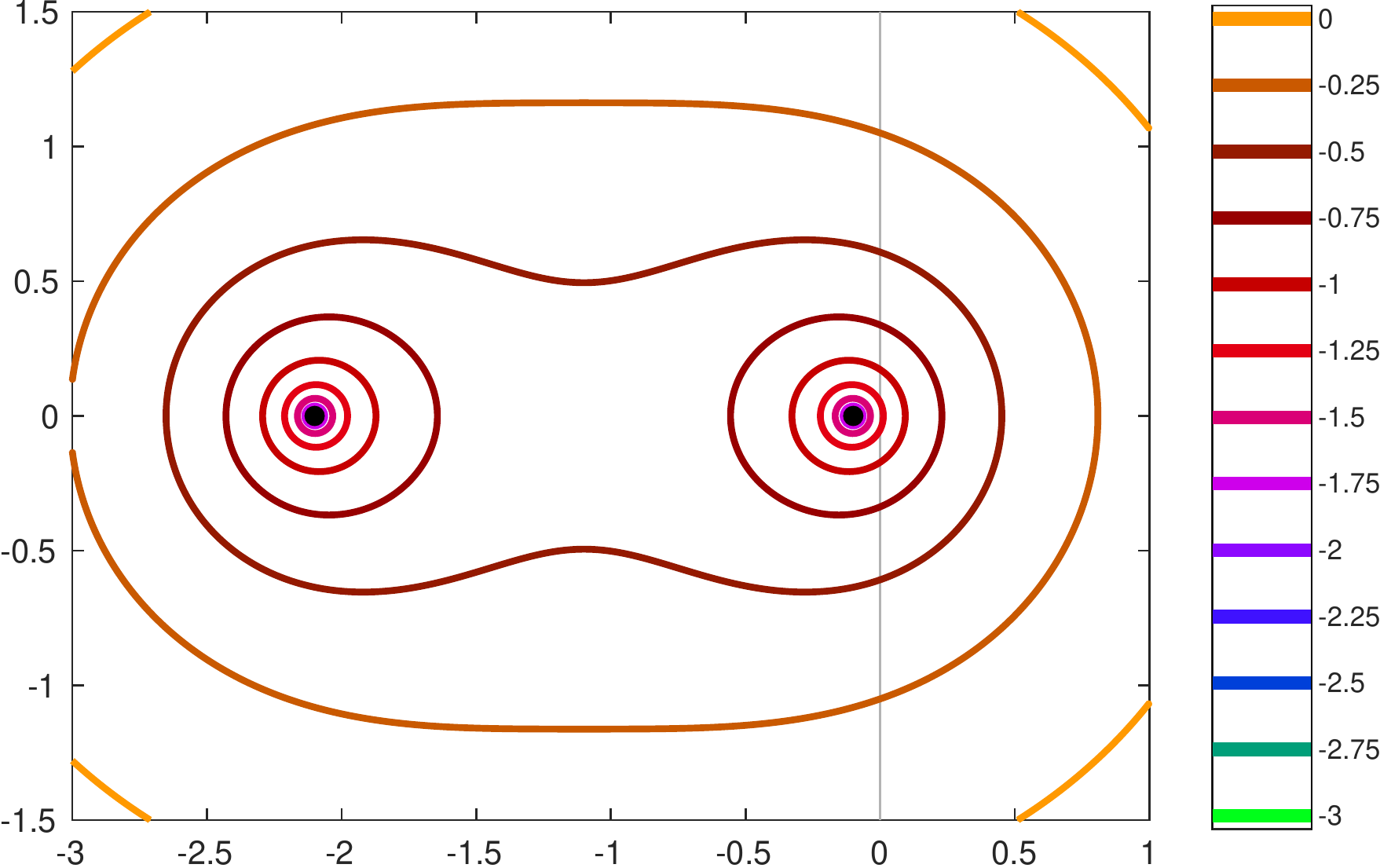} \quad     \includegraphics[width=2.15in]{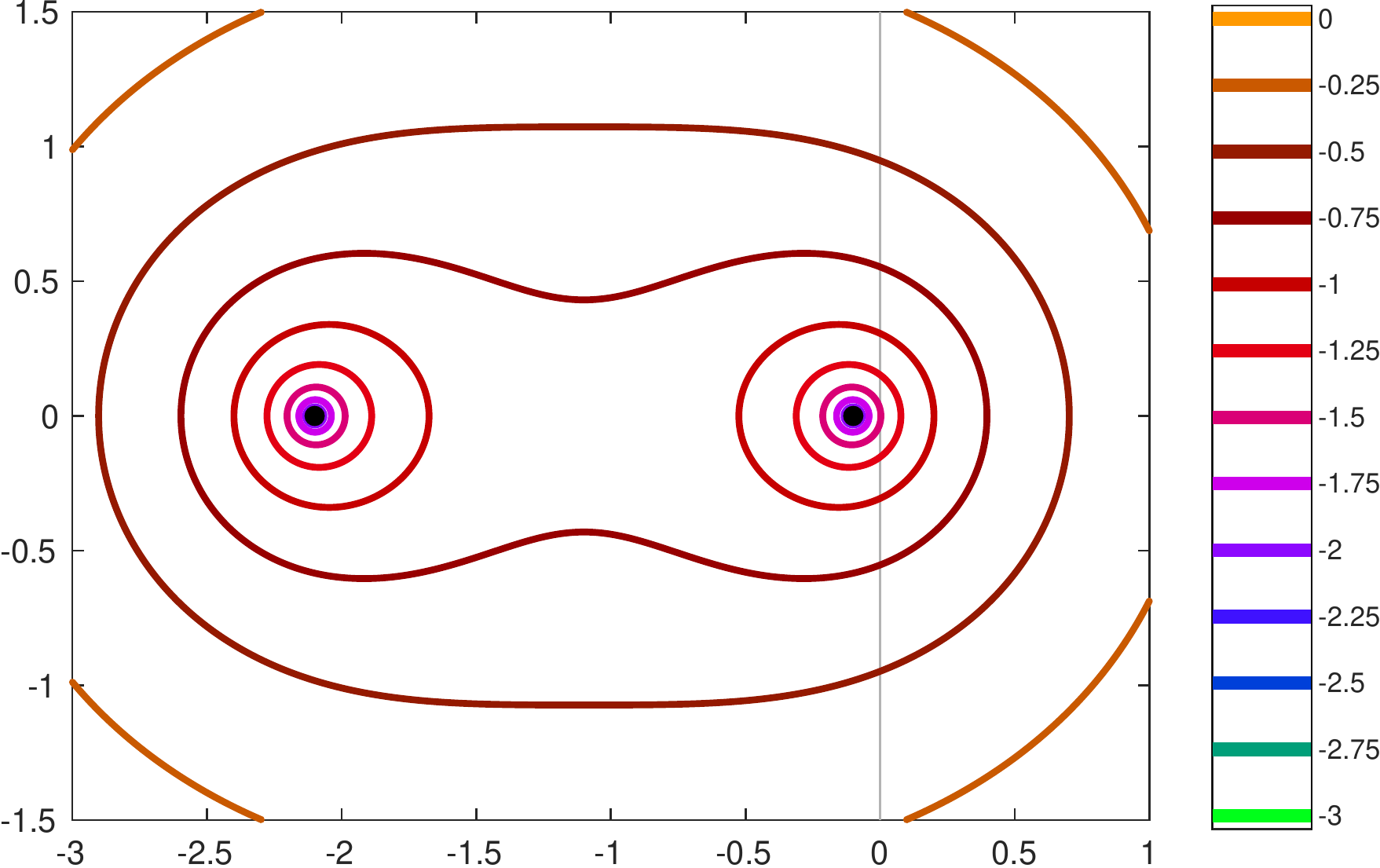}
    \begin{picture}(0,0)
       \put(-214,10){(a)}
       \put(-39,10){(b)}
    \end{picture}

    \vspace*{7pt}
    \includegraphics[width=2.15in]{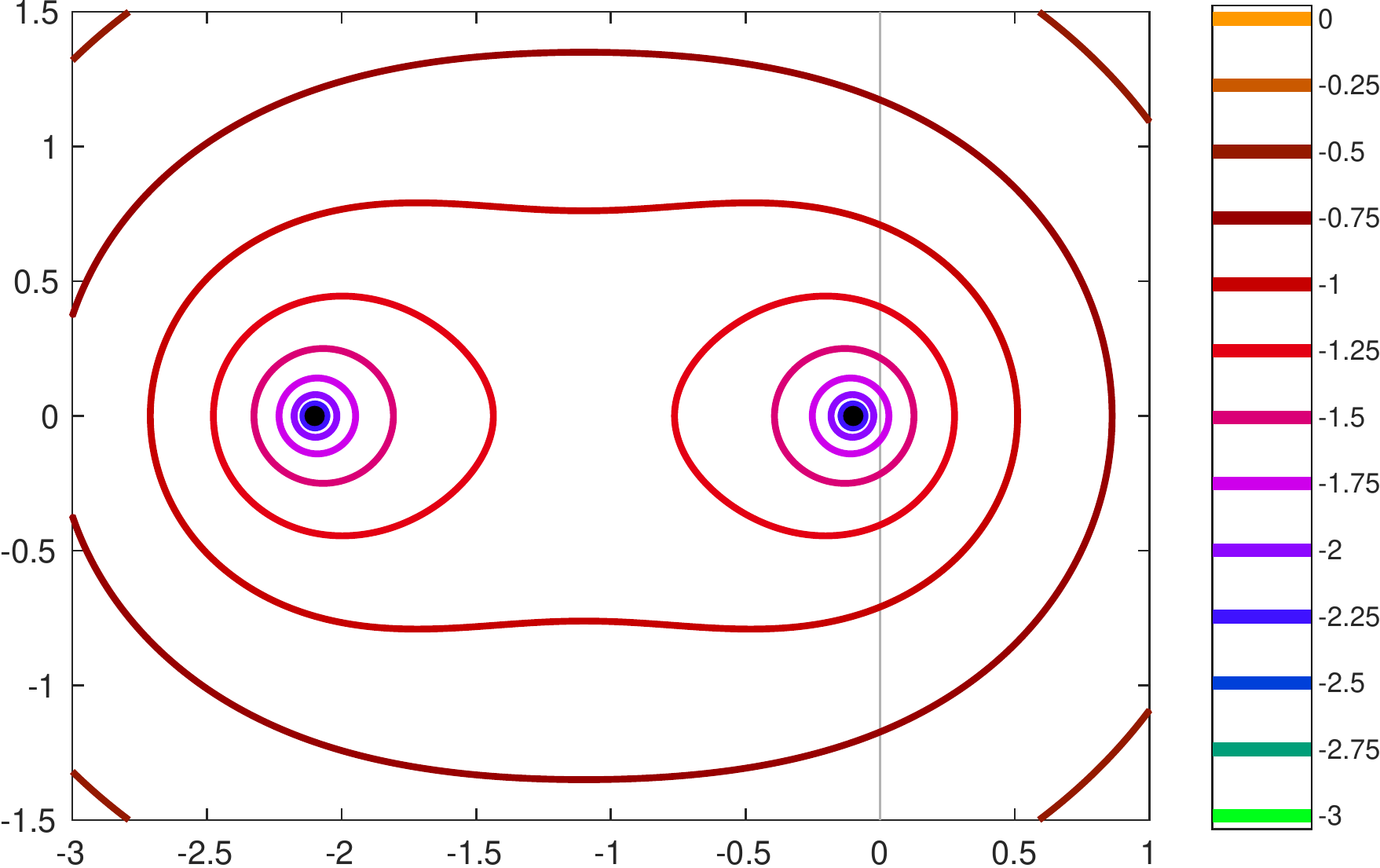} \quad     \includegraphics[width=2.15in]{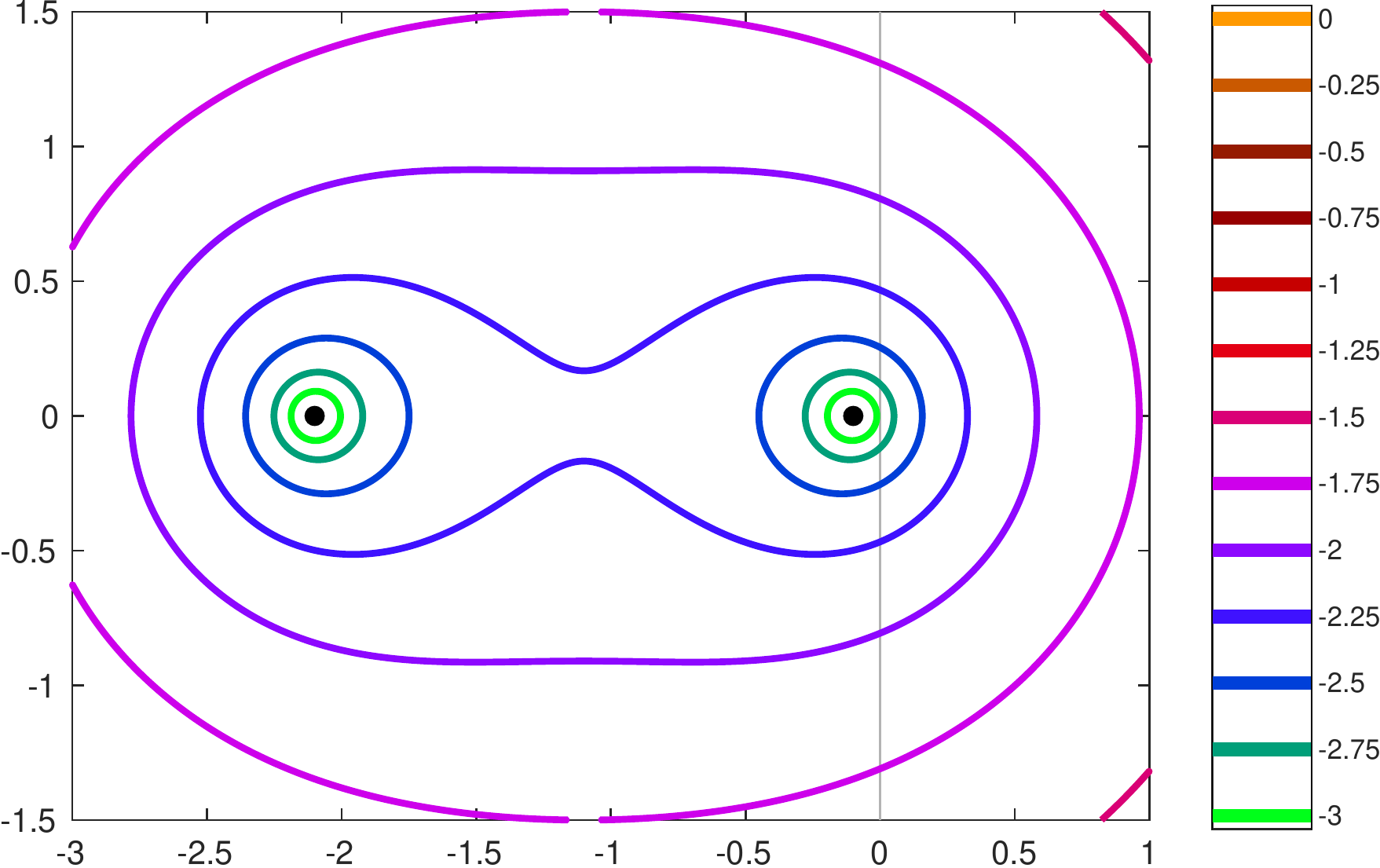}
    \begin{picture}(0,0)
       \put(-214,10){(c)}
       \put(-42,10){(d)}
    \end{picture}
\end{center}

\vspace*{-5pt}
    \caption{\label{fig:twodim1b}
    The pseudospectra $\psa(\BA)$ (top), compared to 
    $\psa(\LL^{-1}\sLL)$ for four Loewner realizations of the 
    system~\eqref{eq:twodim1} using the interpolation points in Table~\ref{tbl:twodim1}.
    In all cases $\sigma(\LL^{-1}\sLL) = \sigma(\BA)$, but the pseudospectra of
    $\LL^{-1}\sLL$ are all quite a bit larger than $\psa(\BA)$.}
\end{figure}

\begin{figure}[p]
\begin{center}
\includegraphics[width=3in]{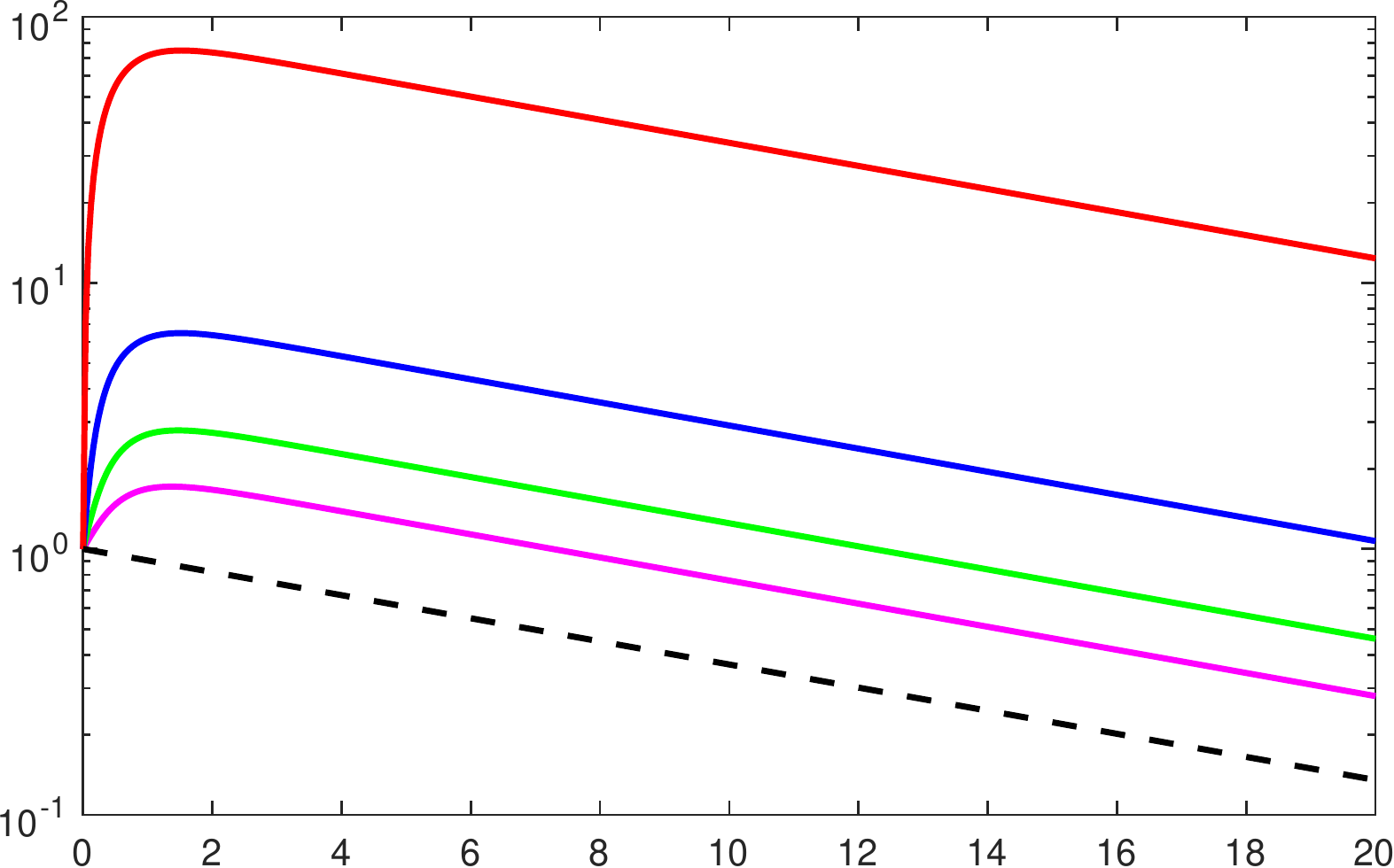}
\begin{picture}(0,0)
\put(-55,-8){$t$}
\put(-250,75){$\|\eop^{t@\LL^{-1}\sLL}\|$}
\put(-92,17){$\|\eop^{t\BA}\|$}
\put(-88,33.5){(a)}
\put(-88,52.5){(b)}
\put(-88,68){(c)}
\put(-88,100){(d)}
\end{picture}
\end{center}

\vspace*{-5pt}
\caption{\label{fig:twodim1_exp}
Evolution of the norm of the solution operator for the original system (black dashed line) and the four interpolating Loewner models.  The instability revealed by the pseudospectra in Figure~\ref{fig:twodim1b} correspond to transient growth in the Loewner systems.}
\end{figure}
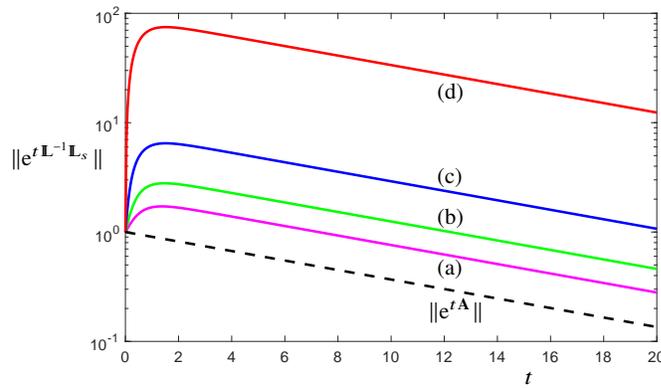

Contrast these results with the standard pseudospectra of $\BA$ itself, 
$\psa(\BA) = \psa^{(1,0)}(\BA,\BI)$ shown at the top of Figure~\ref{fig:twodim1b}.
Since $\BA$ is real symmetric (hence normal), $\psa(\BA)$ is the union
of open $\eps$-balls surrounding the eigenvalues.  
Figure~\ref{fig:twodim1b} compares these pseudospectra to 
$\psa(\LL^{-1}\sLL)$, which give insight about the transient behavior of solutions to 
$\LL \dot\Bx(t) = \sLL\Bx(t)$, e.g., via the bound~(\ref{eq:kreiss}).
The top plot shows $\psa(\BA)$, whose rightmost extent in the complex plane is
always $\eps-0.1$: no transient growth is possible for this system.
However, in all four of the Loewner realizations, $\psa(\LL^{-1}\sLL)$ extends
more than $\eps$ into the right-half plane for $\eps=10^0$ (orange level curve),
indicating by~(\ref{eq:kreiss}) that transient growth must occur
for some initial condition.
Figure~\ref{fig:twodim1_exp} shows this growth for all four realizations:
\emph{the more remote interpolation points lead to Loewner realizations with
greater transient growth}.

\section{Example 2: partitioning interpolation points and noisy data} \label{sec:fullrank2}

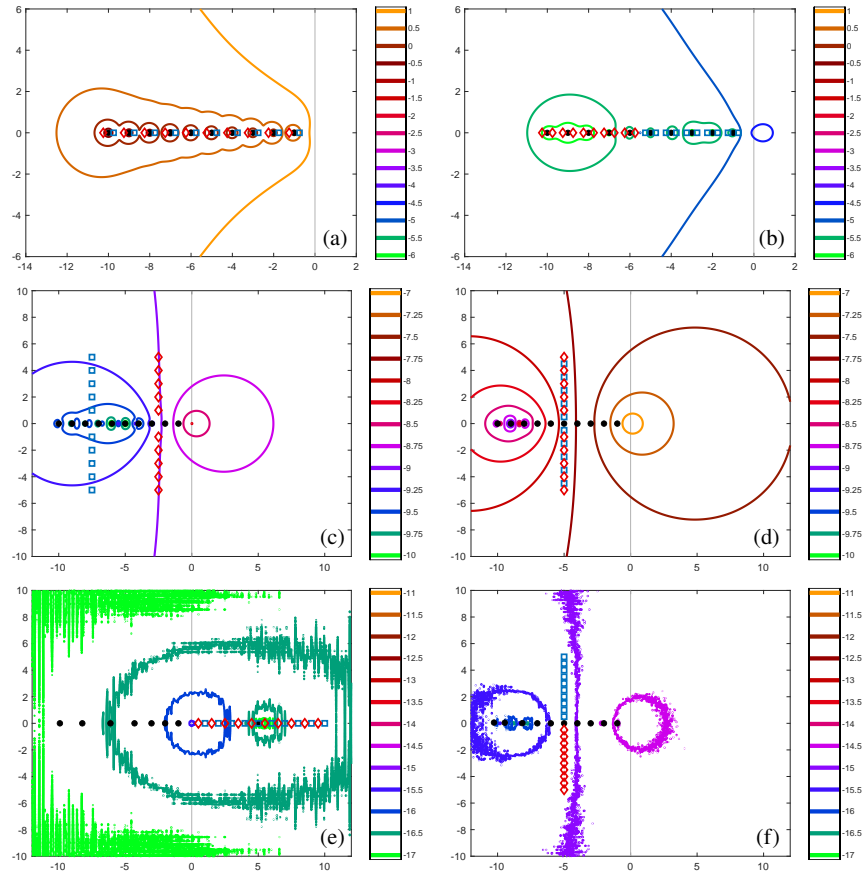
\begin{figure}[t!]
\begin{center}
    \includegraphics[width=2.15in]{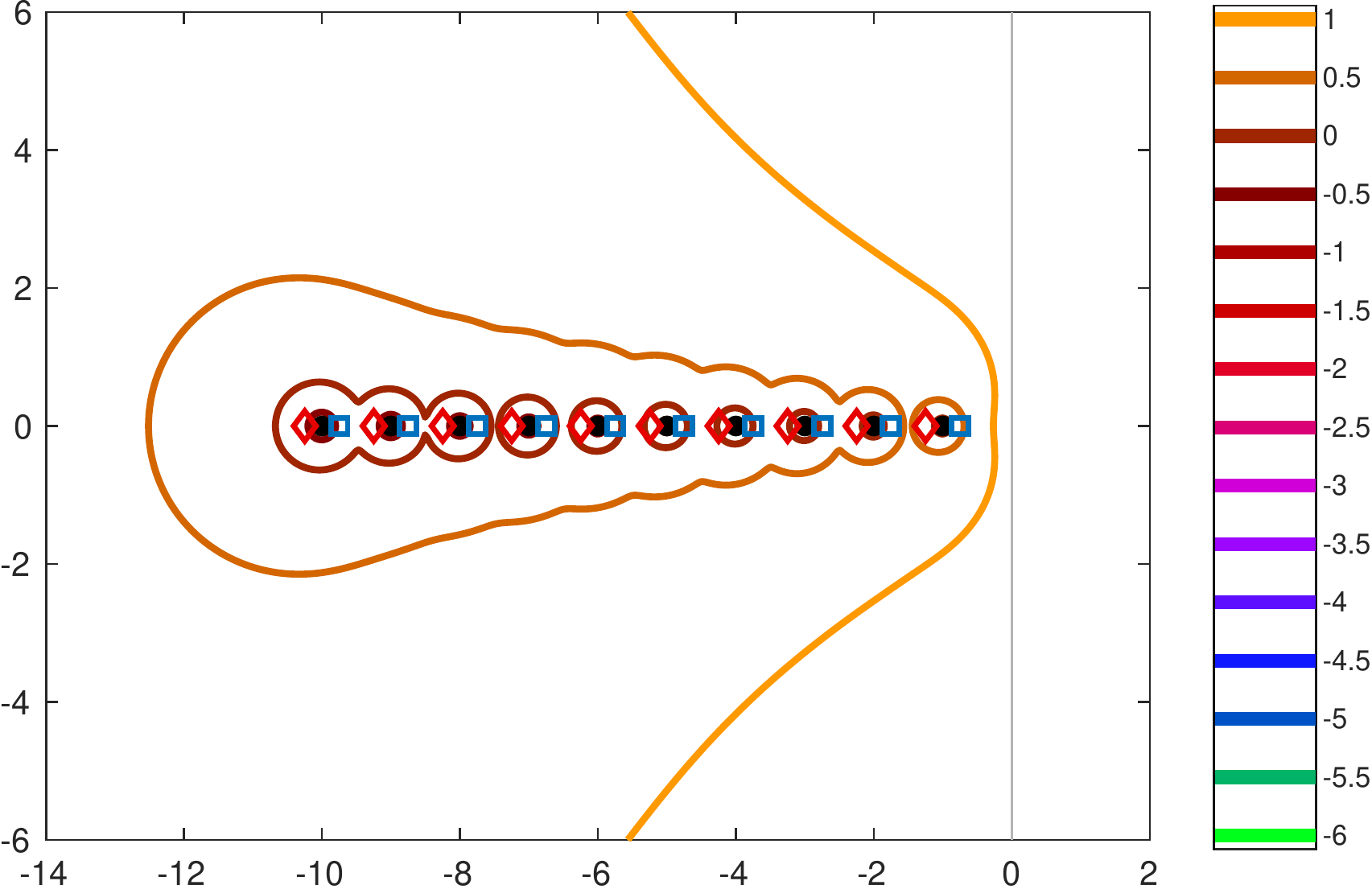} \quad     \includegraphics[width=2.15in]{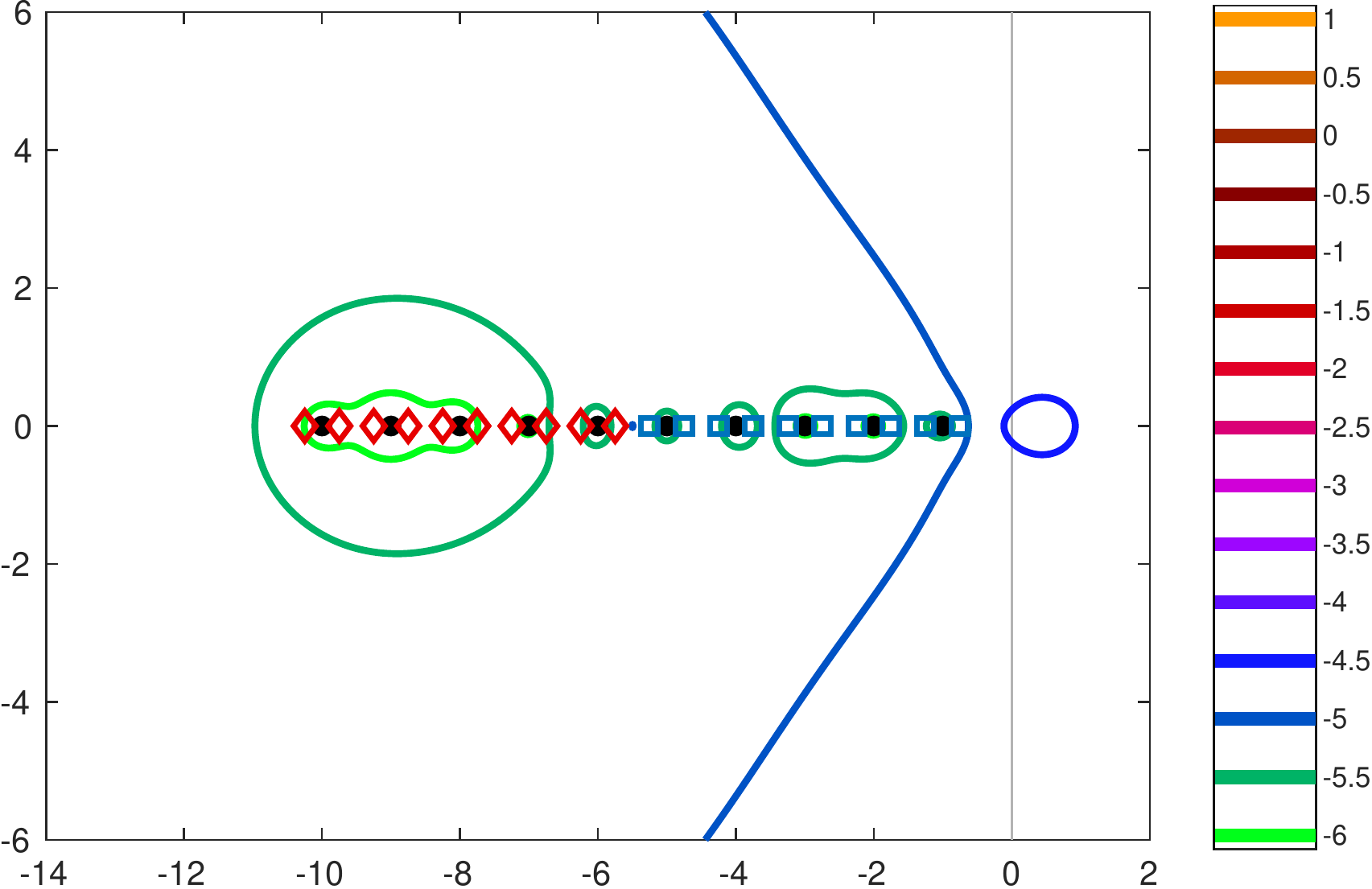}
    \begin{picture}(0,0)
       \put(-206,10){(a)}
       \put(-41,10){(b)}
    \end{picture}

    \vspace*{5pt}
    \includegraphics[width=2.15in]{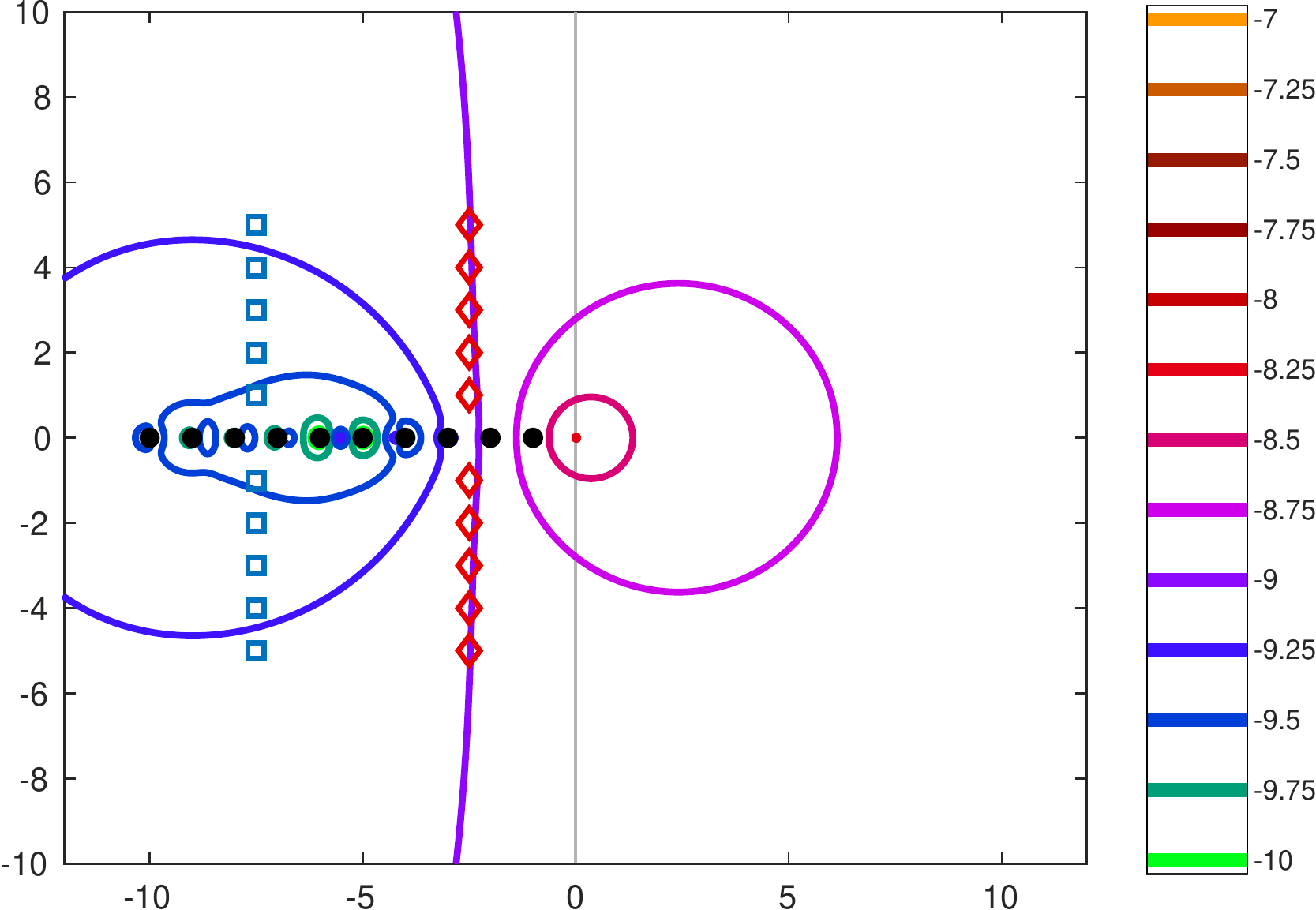} \quad     \includegraphics[width=2.15in]{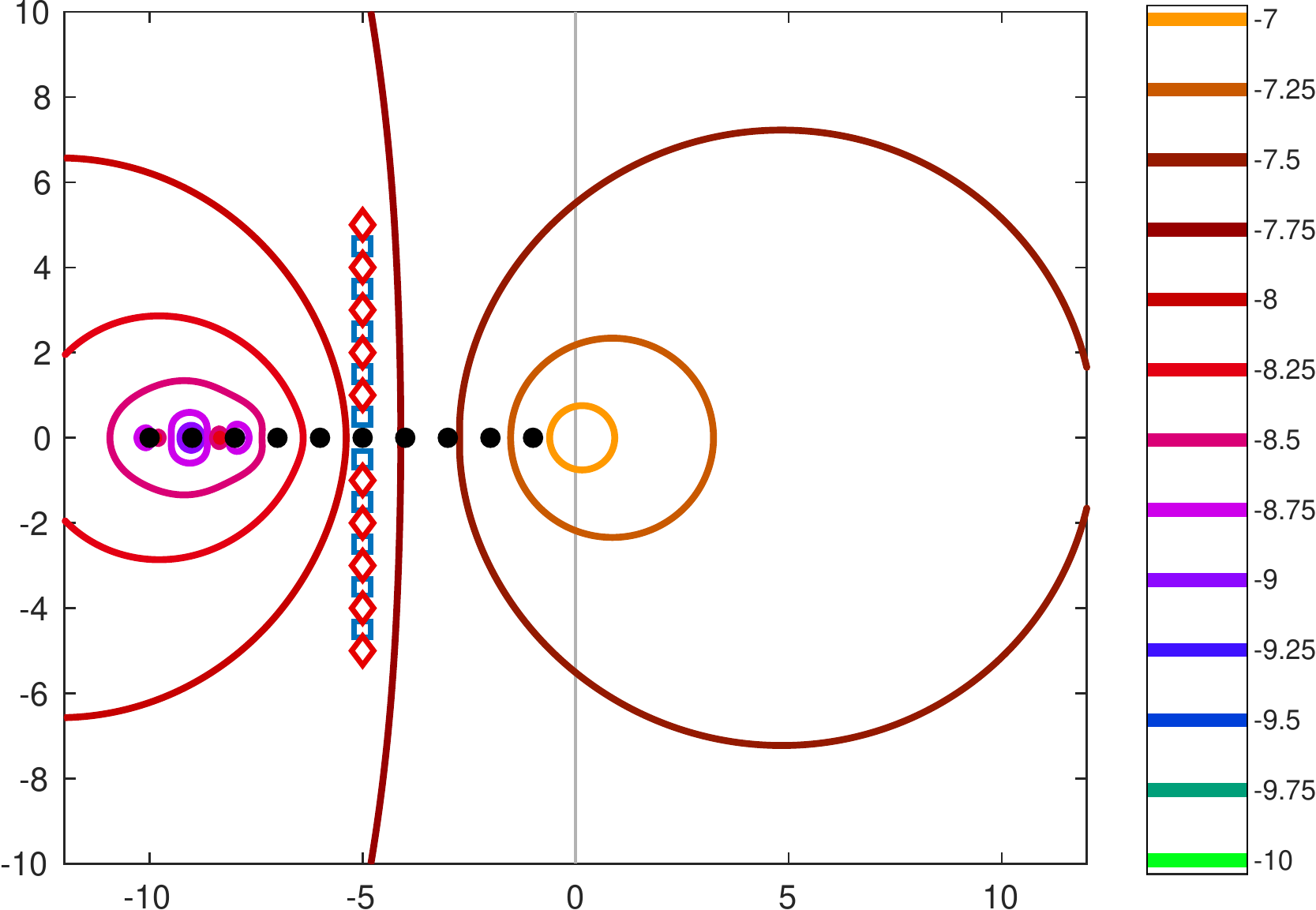}
    \begin{picture}(0,0)
       \put(-207,10){(c)}
       \put(-43,10){(d)}
    \end{picture}

    \vspace*{5pt}
    \includegraphics[width=2.15in]{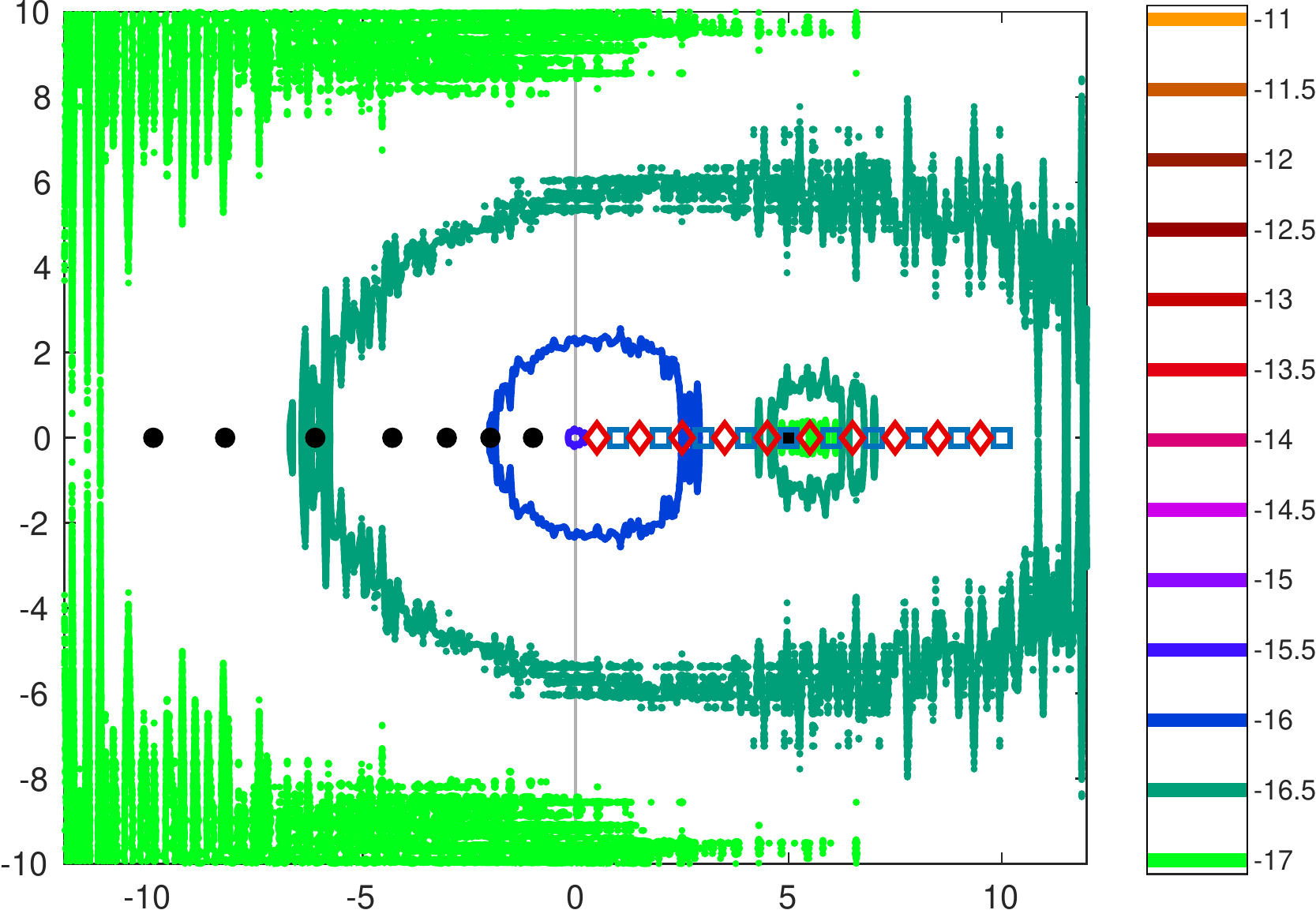} \quad     \includegraphics[width=2.15in]{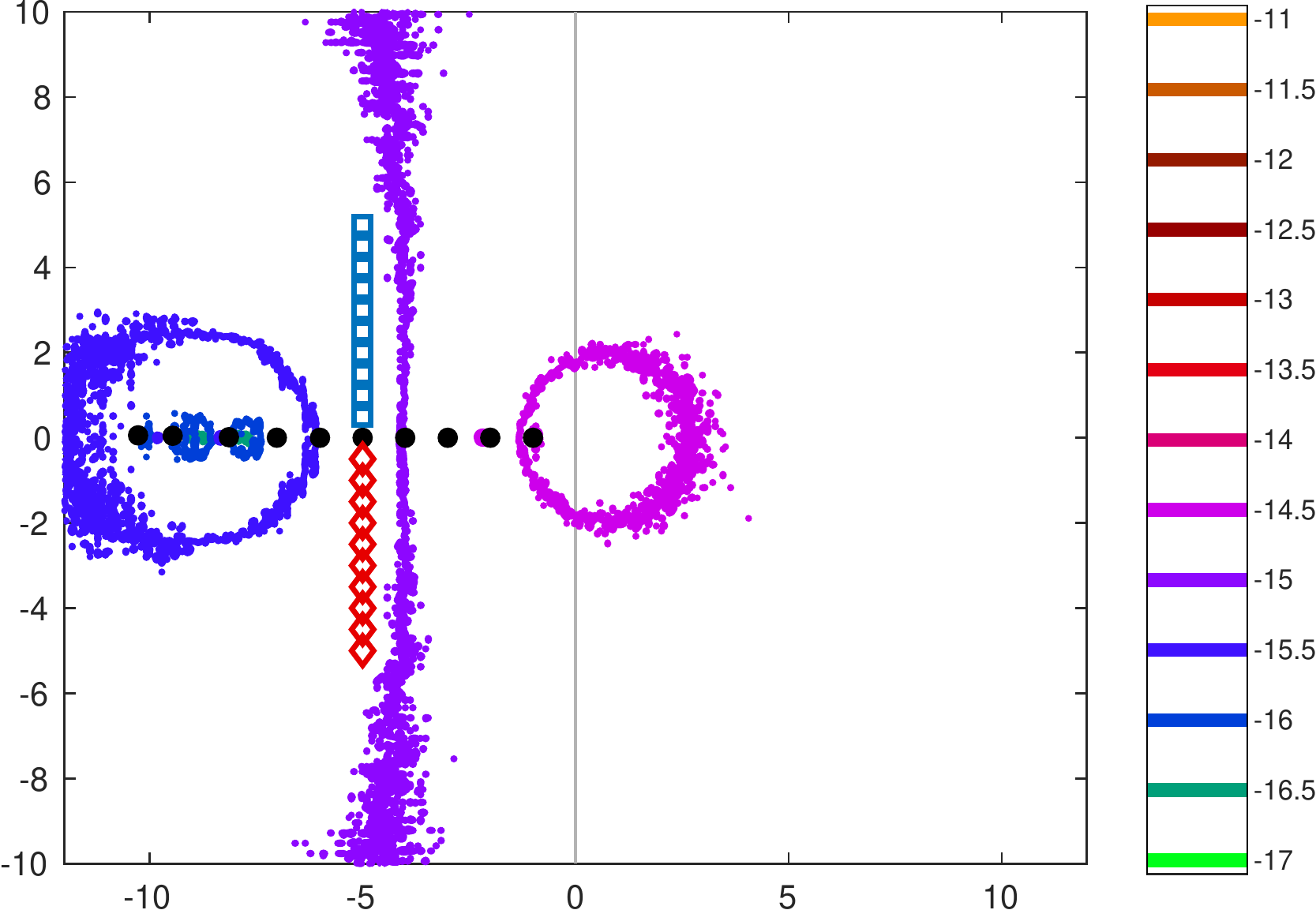}
    \begin{picture}(0,0)
       \put(-207,10){(e)}
       \put(-42,10){(f)}
    \end{picture}
\end{center}
    \caption{\label{fig:siso_rank}
    Pseudospectra $\psa^{(1,1)}(\sLL,\LL)$ for six Loewner realizations of a SISO system of order $n=10$
    with poles $\spec(\BA) = \{-1,-2,\ldots, -10\}$.  
Black dots show computed eigenvalues of the 
pencil $z@\LL-\sLL$ (which should agree with $\spec(\BA)$, but a few are off axis in plot~(e)); 
blue squares show the right interpolation points $\{\lambda_i\}$;
red diamonds show the left points $\{\mu_j\}$.}
\end{figure}

To further investigate how the interpolation points influence eigenvalue stability,
for the Loewner pencil, consider the SISO system of order~10 given by
\[ \BA = {\rm diag}(-1, -2, \ldots, -10), \quad
   \BB = [1, 1, \ldots, 1]^T, \quad
   \BC = [1, 1, \ldots, 1].\]
Figure~\ref{fig:siso_rank} shows $\psa^{(1,1)}(\sLL,\LL)$ for six configurations of the
interpolation points.  Plots~(a) and~(b) use the points 
$\{-10.25, -9.75, -9.25, \ldots, -1.25, -0.75\}$; 
in plot~(a) the left and right points interlace, suggesting slower decay of the 
singular values of $\LL$, as discussed in Section~\ref{sec:syl};
in plot~(b) the left and right points are separated, leading to faster
decay of the singular values of $\LL$ and considerably larger pseudospectra.
(The Beckermann--Townsend bound~\cite[Cor.~4.2]{BT17} would apply to this case.)
Plot~(c) further separates the left and right points, giving even larger pseudospectra.
In plots~(d) and (f), complex interpolation points $\{-5\pm0.5@\iop, -5\pm1.0@\iop,\ldots,-5\pm5@\iop\}$
are interleaved (while keeping conjugate pairs together)~(d) and separated~(f): 
the latter significantly enlarges the pseudospectra.
Plot~(f) uses the same relative arrangement that gave such nice results in plot~(a)
(the singular value bound~\eqref{eq:syl} is the same for (a) and (e)), but their locations relative
to the poles of the original system differ.  The pseudospectra are now much larger, showing 
that a large upper bound in~(\ref{eq:syl}) is not alone enough to guarantee small 
pseudospectra.  (Indeed, the pseudospectra are so large in~(e) and~(f) that the plots are
dominated by numerical artifacts of computing $\|(z\BE-\BA)^{-1}\|$.)
\emph{Pseudospectra reveal the great influence interpolation point location and partition
can have on the stability of the realized pencils.}

\begin{figure}[t!]

\begin{center}
    \includegraphics[width=2in]{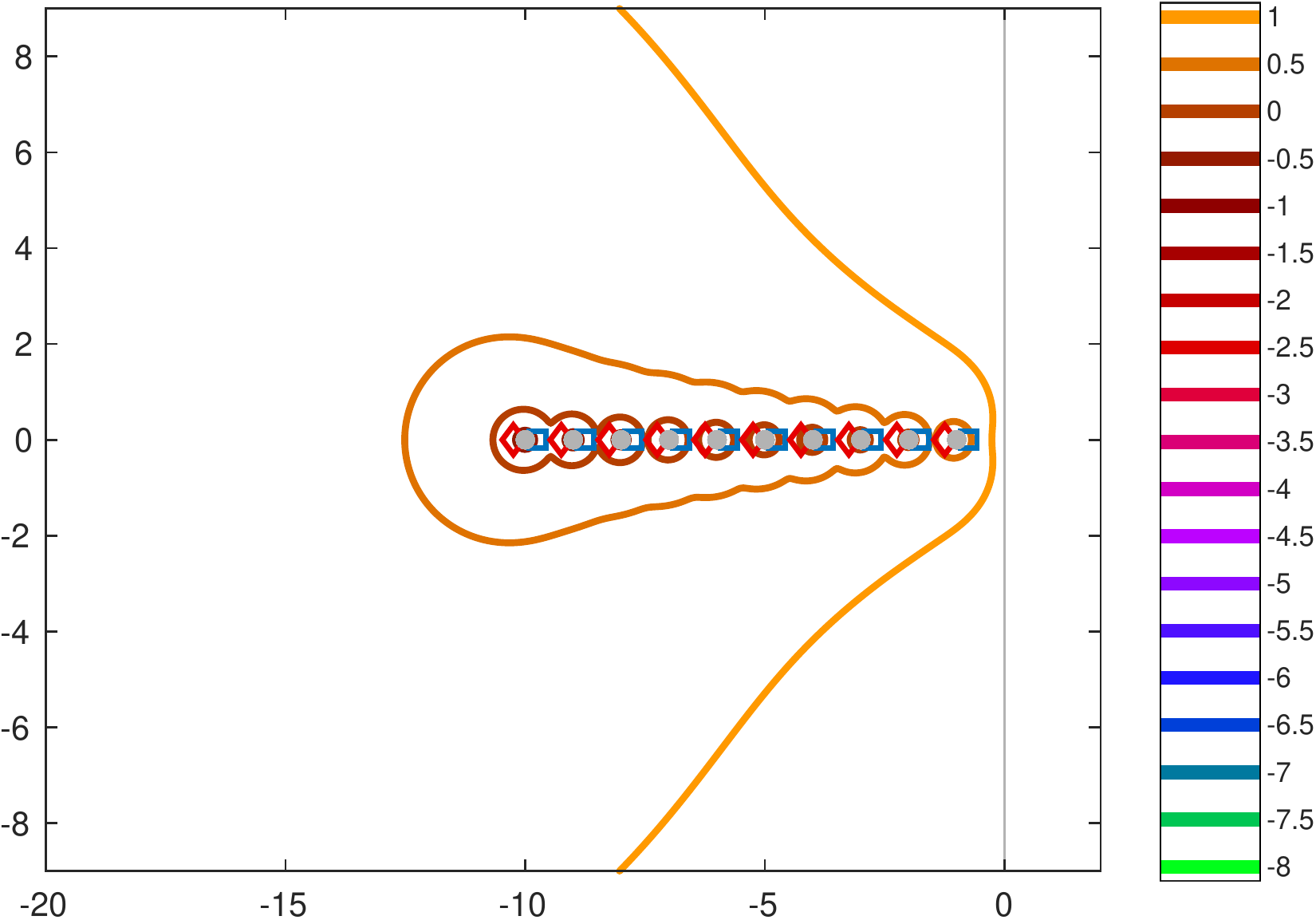}
    \begin{picture}(0,0) 
       \put(-138,10){\tiny noise = $10^{-1}$}
    \end{picture}
    \qquad
    \includegraphics[width=2in]{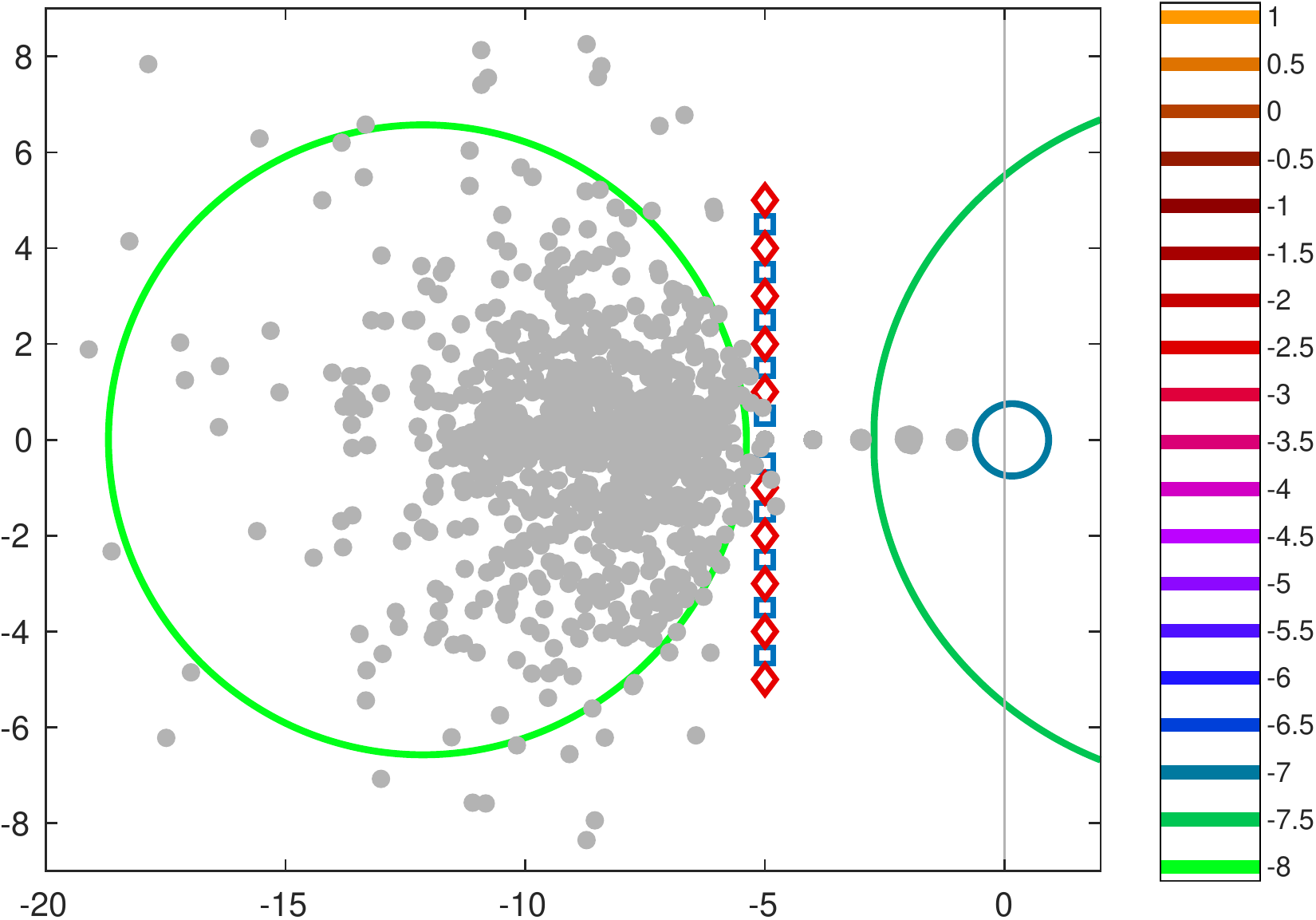} 
    \begin{picture}(0,0) 
       \put(-138,8){\tiny  noise = $5\cdot10^{-9}$}
    \end{picture}
\end{center}

\vspace*{-5pt}
\caption{\label{fig:siso_pert}
Eigenvalues of the perturbed Loewner pencil $\wh{\LL}_s-z\wh{\LL}$ (gray dots),
constructed from measurements that have been perturbed by random complex noise of 
magnitude $10^{-1}$ (left) and $5\cdot10^{-9}$ (right) (1,000 trials).  
As the pseudospectra in Figure~\ref{fig:siso_rank}(a,d) indicate, 
the interleaved interpolation points on the left are remarkably stable, while
the similarly interleaved complex interpolation points on the right give 
a Loewner pencil that is incredibly sensitive to small changes in the data.
}
\end{figure}
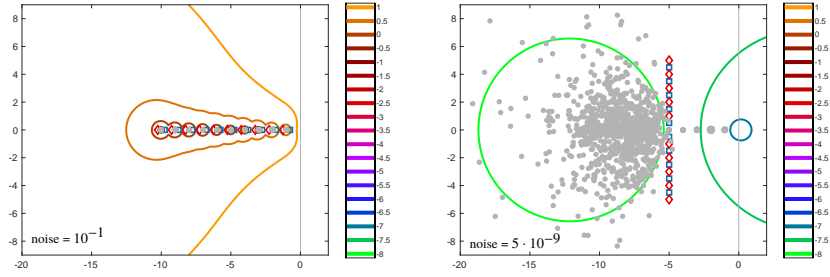

Pseudospectra also give insight into the consequences of inexact measurement data.
Consider the following experiment.  Take the scenario in Figure~\ref{fig:siso_rank}(a),
the most robust of these examples.
Subject each right and left measurement $\{\Bw_1,\ldots, \Bw_\rho\} \subset \C$ and
$\{\Bv_1,\ldots, \Bv_\nu\} \subset \C$ to random complex noise of magnitude $10^{-1}$,
then build the Loewner pencil $\wh{\LL}_s-z\wh{\LL}$ from this noisy data.
How do the badly polluted measurements affect the computed eigenvalues?
Figure~\ref{fig:siso_pert} shows the results of 1,000~random trials,
which can depart from the true matrices significantly:
\[ 3.06 \le \|\wh{\LL}_s - \sLL\| \le 5.88, \qquad 0.49 \le \|\wh{\LL}-\LL\| \le 0.63. \]
Despite these large perturbations, the recovered eigenvalues are remarkably
accurate:  in $99.99\%$ of cases, the eigenvalues have absolute accuracy of at least $10^{-2}$, 
indeed more accurate than the measurements themselves.  
The pseudospectra in Figure~\ref{fig:siso_rank}(a) suggest good robustness
(though the pseudospectral level curves are pessimistic by one or two orders of magnitude).
Contrast this with the complex interleaved interpolation points used in Figure~\ref{fig:siso_rank}(d).
Now we only perturb the data by a small amount, $5\cdot 10^{-9}$, for which the perturbed Loewner
matrices (over 1,000 trials) satisfy
\[ 1.52\cdot 10^{-7} \le \|\wh{\LL}_s - \sLL\| \le 2.03\cdot 10^{-7}, 
   \qquad 2.57\cdot 10^{-8} \le \|\wh{\LL}-\LL\| \le 3.14\cdot 10^{-8}. \]
With this mere hint of noise, the eigenvalues of the recovered system erupt:
only $36.73\%$ of the eigenvalues are correct to two digits.
(Curiously, $-4$ and $-5$ are always computed correctly, while $-8$, $-9$, and $-10$
are never computed correctly.)  
The pseudospectra indicate that the leftmost eigenvalues
are more sensitive, and again hint at the effect of the perturbation (though off
by roughly an order of magnitude in $\eps$).%
\footnote{One could inspect the $\eps=1$ level curve of $\sigma_\eps^{(\gamma_\star, \delta_\star)}(\sLL,\LL)$
for $\gamma_\star = \max \{ \|\sLL-\wh{\LL}_s\|\}$ and $\delta_\star = \max \{ \|\LL-\wh{\LL}\|\}$.}
Measurements of real systems (or even numerical simulations of nontrivial systems)
are unlikely to produce such high accuracy; pseudospectra
can reveal the virtue or folly of a given interpolation point configuration.

\begin{figure}[t!]
\begin{center}
    \includegraphics[width=2in]{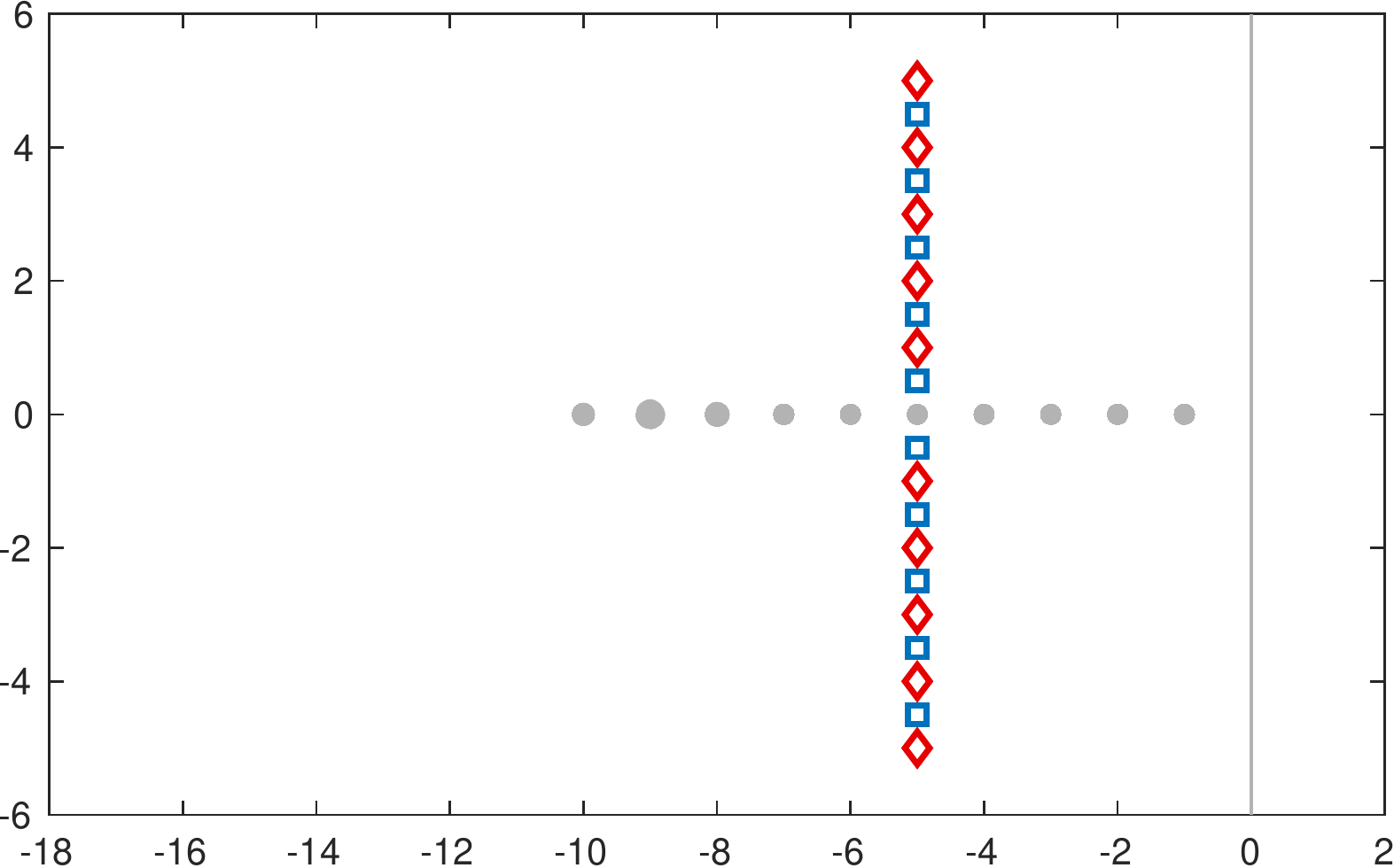} \qquad
    \includegraphics[width=2in]{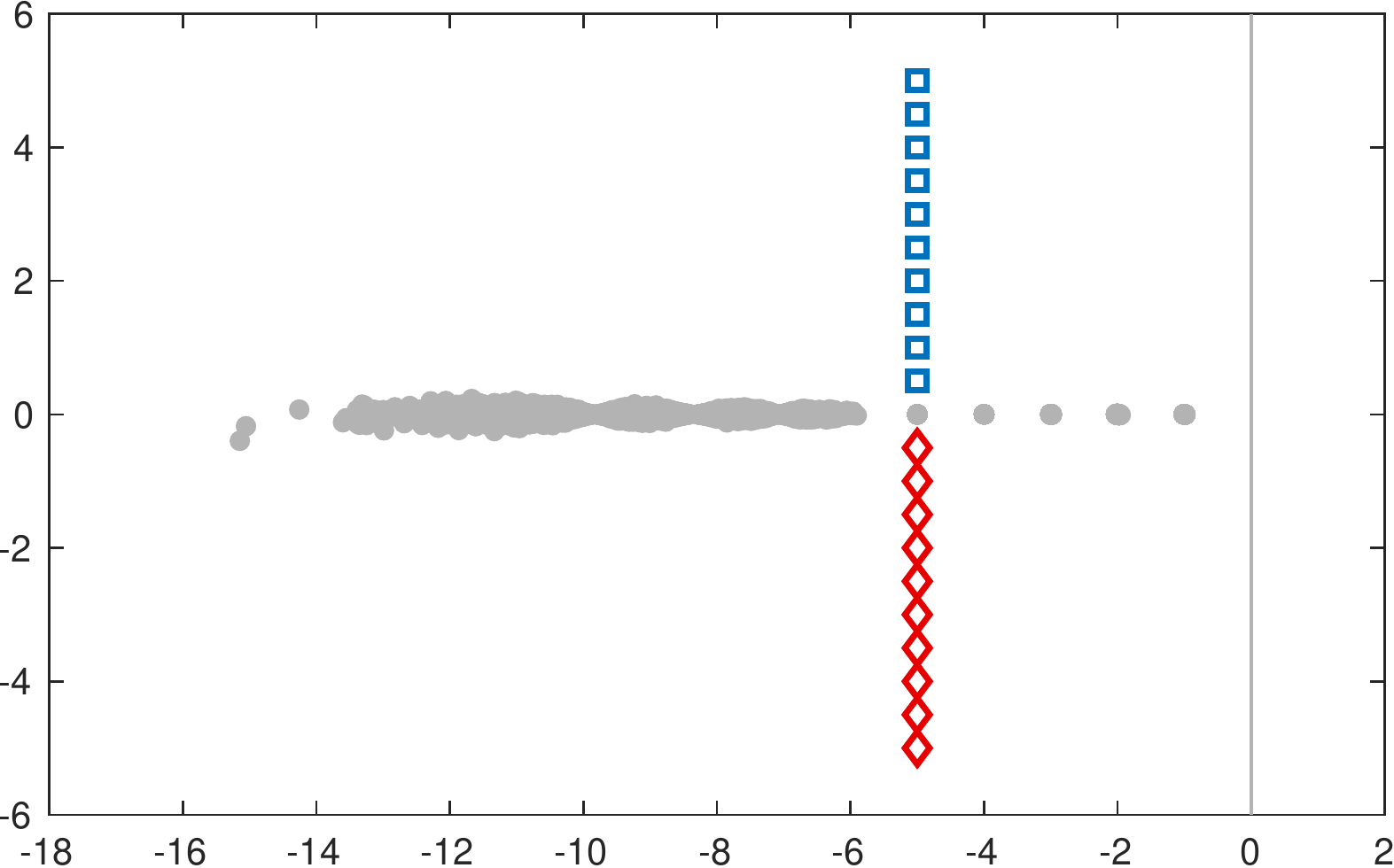} 
\end{center}
\vspace*{-5pt}
\caption{\label{fig:siso_rank6b_pert}
Eigenvalues of the perturbed Loewner pencil $\wh{\LL}_s-z\wh{\LL}$ (gray dots),
constructed from measurements that have been perturbed by random complex noise of 
magnitude $10^{-10}$ (10,000 trials).  As suggested by the pseudospectra plots,
the interleaved interpolation points (left) are more robust to perturbations 
than the separated points (right), though the difference is not as acute as 
suggested by Figure~\ref{fig:siso_rank}(d,f).
For example, in these 10,000 trials, the least stable pole ($-9$) is computed 
accurately (absolute error less than $.01$) in $10.97\%$ of trials on the left, 
and $0.29\%$ on the right.
}
\vspace*{-9pt}
\end{figure}

In these simple experiments, pseudospectra have been most helpful for indicating
the sensitivity of eigenvalues when the left and right interpolation points are 
favorably partitioned (e.g., interleaved).
They seem to be less precise at predicting the sensitivity to noise of poor
left/right partitions of the interpolation points. 
Figure~\ref{fig:siso_rank6b_pert} gives an example, based on the two partitions
of the same interpolation points in Figure~\ref{fig:siso_rank}(d,f).
The pseudospectra suggest that the eigenvalues for plot~(f) should be much more
sensitive to noise than those for the interleaved points in plot~(d).
In fact, the configuration in plot~(f) appears to be only marginally less stable
to noise of size $10^{-10}$, over 10,000 trials.
This is a case where one could potentially glean additional insight from 
\emph{structured} Loewner pseudospectra.
\vspace*{-10pt}

%

\vspace*{-9pt}

\section{Conclusion}

\vspace*{-9pt}
Pseudospectra provide a tool for analyzing the stability
of eigenvalues of Loewner matrix pencils.
Elementary examples show how pseudospectra can inform the
selection and partition of interpolation points, and bound the eigenvalues of Loewner
pencils in the presence of noisy data.  Using a different approach to pseudospectra,
we showed that while the realized Loewner pencil matches the poles of the
original system, it need not replicate transient dynamics of $\dot\Bx(t) = \BA\Bx(t)$;
pseudospectra can reveal potential transient growth, which
varies with the interpolation points.

In this initial study we have intentionally used simple examples involving small,
symmetric $\BA$.  
Realistic examples, e.g., with complex poles, nonnormal $\BA$, multiple inputs and outputs,
and rank-deficient Loewner matrices, will add additional complexity.
Moreover, we have only used the Mayo--Antoulas interpolation theory to realize a 
system whose order is known; we have not addressed pseudospectra of the reduced 
pencils~(\ref{eq:redLoew}) in the context of data-driven model reduction.

\emph{Structured Loewner pseudospectra} provide another avenue for future study.  
Structured \emph{matrix pencil} pseudospectra 
have not been much investigated, especially with Loewner structure.  
Rump's results for standard pseudospectra~\cite{Rum06} suggest the following problem; its positive resolution would imply that
the Loewner pseudospectrum $\psagd(\sLL,\LL)$ matches the structured Loewner matrix
pseudospectrum.

Given any $\eps, \gamma, \delta > 0$, Loewner matrix $\LL$ and associated shifted
Loewner matrix $\sLL$, suppose $z \in \psagd(\sLL,\LL)$.
\emph{Does there exist some Loewner matrix $\wh{\LL}$ and associated shifted Loewner
matrix $\wh{\LL}_s$ such that 
\[ \|\wh{\LL}_s - \sLL\| < \epsilon@\gamma, \qquad 
   \|\wh{\LL} - \LL\| < \epsilon@\delta\]
and 
$z\in \sigma(\wh{\LL}_s,\wh{\LL})$}\,?


\newpage
\begin{acknowledgement}
This work was motivated by a question posed by Thanos Antoulas,
who we thank not only for suggesting this investigation, but also for
his many profound contributions to systems theory and his inspiring 
teaching and mentorship.  We also thank Serkan Gugercin for helpful comments.
(Mark Embree was supported by the U.S. National Science Foundation 
under grant DMS-1720257.)
\end{acknowledgement}

\section*{Appendix}
We provide a MATLAB implementation that computes the inverse iteration vectors $\Bu$
in (\ref{eq:inviter}) in only $\CO(n^2)$ operations
by exploiting the Cauchy-like rank displacement structure of the Loewner pencil,
as shown in~(\ref{eq:inviter2}). Namely, we start from the general $\CO(n^3)$ MATLAB code 
from~\cite[p.~373]{TE05} and modify it to account for the Loewner structure,
and hence achieve $\CO(n^2)$ efficiency.
This code computes $\|(z\LL-\sLL)^{-1}\|$ for a fixed $z$.  To compute pseudospectra,
one applies this algorithm on a grid of $z$ values.  In that case, the 
$\CO((m+p)n^2)$ structured LU factorization in the first line need only be computed \emph{once} 
for all $z$ values (just as the standard algorithm computes an $\CO(n^3)$ simultaneous 
triangularization using the QZ algorithm once for all $z$).

{\small
\begin{verbatim}
  [L1,U1,piv1] = LUdispPiv(mu,lambda,[V' L'],[R.' -W.']);
  L1t = L1'; U1t = U1';
  z = 1./(z-lambda); Upz = z.*(U1\(L1\V(:,piv1)'));

  [L2,U2,piv2] = lu(eye(m)-R*Upz,'vector'); R2 = R(piv2,:);
  L2t = L2'; U2t = U2'; R2t = R2'; Upzt = Upz';
  applyTheta = @(x) x+Upz*(U2\(L2\(R2*x)));
  applyThetaTranspose = @(x) x+R2t*(L2t\(U2t\(Upzt*x)));

  sigold = 0;
  for it = 1:maxit
      u = U1\(L1\u(piv1));
      u = conj(z).*applyThetaTranspose(applyTheta(z.*u));
      u(piv1) = L1t\(U1t\u);

      sig = 1/norm(u);
      if abs(sigold/sig-1) < 1e-2, break, end
      u = sig*u; sigold = sig;
  end
  sigmin = sqrt(sig);
\end{verbatim}}
The function {\tt LUdispPiv} computes the LU factorization (with partial pivoting)
of the Loewner matrix $\LL$ in $\CO((m+p)n^2)$ operations.
The Loewner matrix is not formed explicitly; instead, the function uses the raw interpolation data
$\lambda_i, \Br_i, \Bw_i$ and $\mu_j, \Bell_j, \Bv_j$.
The implementation details for {\tt LUdispPiv} can be found in \cite[sect.~12.1]{GV12}.
The LU factorization of $\BI - \BR\BV^*(z) \in \Cmm$ is given by {\tt L2} and {\tt U2}, while the first
three lines of the loop represent the computation of $\LL^{-*}\BTheta(z)^*\BTheta(z)\LL^{-1}\Bu$, 
as defined in~(\ref{eq:inviter2}). Note the careful grouping of terms and the use of elementwise
multiplication {\tt .* } to keep the total operation count at $\CO(n^2)$.


\end{document}